\documentclass[11pt]{article}
\usepackage{amssymb, amsthm, amsmath, amscd}
\newtheorem{thm}{Theorem}[section]
\newtheorem{lem}[thm]{Lemma}
\newtheorem{prop}[thm]{Proposition}
\newtheorem{cor}[thm]{Corollary}
\newtheorem{NN}[thm]{}
\theoremstyle{definition}\newtheorem{df}[thm]{Definition}
\theoremstyle{definition}
\theoremstyle{definition}

\renewcommand{\phi}{\varphi}

\newcommand{\Z}{\mathbb{Z}}
\newcommand{\Q}{\mathbb{Q}}
\newcommand{\R}{\mathbb{R}}
\newcommand{\C}{\mathbb{C}}
\newcommand{\T}{\mathbb{T}}

\newcommand{\morp}{contractive completely positive linear map}

\newcommand{\hm}{homomorphism}
\newcommand{\dt}{\delta}
\newcommand{\ep}{\epsilon}
\newcommand{\andeqn}{\,\,\,{\rm and}\,\,\,}
\newcommand{\rforal}{\,\,\,{\rm for\,\,\,all}\,\,\,}
\newcommand{\CA}{$C^*$-algebra}
\newcommand{\SCA}{$C^*$-subalgebra}

\newcommand{\af}{{\alpha}}
\newcommand{\bt}{{\beta}}

\newcommand{\D}{\mathbb D}
\newcommand{\beq}{\begin{eqnarray}}
\newcommand{\eneq}{\end{eqnarray}}
\newcommand{\tforal}{\,\,\,\text{for\,\,\,all}\,\,\,}
\newcommand{\tand}{\,\,\,\text{and}\,\,\,}
\title{ AF-embedding of the crossed products of AH-algebras by finitely generated abelian groups }
\author{Huaxin Lin
 }
\date{}

\begin{document}

\maketitle

\begin{abstract}

Let $X$ be a compact metric space and let $\Lambda$ be a $\Z^k$
($k\ge 1$) action on $X.$ We give a solution to a version of
Voiculescu's problem of AF-embedding: The crossed  product
$C(X)\rtimes_{\Lambda}\Z^k$ can be embedded into a unital simple
AF-algebra if and only if $X$ admits a strictly positive
$\Lambda$-invariant Borel probability  measure.

Let $C$ be a unital AH-algebra, let $G$ be a finitely generated abelian group and let
$\Lambda: G\to Aut(C)$ be a monomorphism. We show that  $C\rtimes_{\Lambda} G$ can be embedded into a unital simple AF-algebra if and only if $C$ admits a faithful $\Lambda$-invariant tracial state.

\end{abstract}

\section{Introduction}
Quasidiagonality for crossed product \CA s were studied by Dan
Voiculescu (\cite{V2}), \cite{V3} and \cite{V4}). Quasidiagonality
in \CA s has been studied for a long time in many different point
of views (see \cite{Bn2} for more information, also \cite{DHS},
\cite{Th}, \cite{Hd}, \cite{Ros}, \cite{Sa}, \cite{Z}, \cite{BK1},
\cite{BK2}, \cite{BlD}, \cite{ELP}, \cite{D1} and \cite{D2}, to
name a few). Let $X$ be a compact metric space and $\af$ be a
homeomorphism on $X.$  It was proved by Pimsner (\cite{Pi}) that
$C(X)\rtimes_{\af}\Z$ is quasidiagonal if and only if $\af$ is
pseudo-non-wondering, and if and only if $C(X)\rtimes_{\af}\Z$ can
be embedded into an AF-algebra. A more general  question is when
$A\rtimes_{\af} \Z$ can be embedded into an AF-algebra, where
$\af\in Aut(A)$ is an automorphism on $A$ and where $A$ is a
unital separable \CA\, which can be embedded into a unital simple
AF-algebra.  Let $\Lambda$ be a $\Z^2$ action on $X.$  Dan
Voiculescu (4.6 of \cite{V4} ) asked  when
$C(X)\rtimes_{\Lambda}\Z^2$ can be embedded into an AF-algebra.

For the first question, there are some significant  progresses around the turn of the century.
Let $A$ be an AF-algebra and $\af\in Aut(A)$ be an automorphism.
Nate Brown   (\cite{Bn1}) proved (1998) that $A\rtimes_{\af}\Z$ can be
embedded into an AF-algebra if and only if $A\rtimes_{\af}\Z$ is
quasidiagonal. More importantly, a $K$-theoretical necessary and sufficient
condition for $A\rtimes_{\af}\Z$ being quasidiagonal is also
given there. When $A$ is a unital simple $A\T$-algebra of real rank zero, Matui (around 2001)
(\cite{M1}) showed that $A\rtimes_{\af}\Z$ can always be embedded into a unital simple AF-algebra.
We find that it is particularly important to know when a crossed product can be embedded into a unital simple AF-algebra.
 More recently, it was proved  (\cite{Lnasy}) that,
for any unital AH-algebra $A,$ $A\rtimes_{\af}\Z$ can be embedded into a unital simple AF-algebra if (and only if)
$A$ has a faithful $\af$-invariant tracial state. Note that $A$ is not assumed to have real rank zero, nor simple.

It is much more difficult to answer the Voiculescu question above.
However,  there  were some progresses made in last 15 years or so.
H. Matui (\cite{M1}) showed that under some additional condition on
the action, $C(X)\rtimes_{\Lambda}\Z^2$ can be embedded into an
AF-algebra. N. Brown (\cite{Bn3}) proved that if $A$ is a
UHF-algebra and $\Lambda:\Z^k\to Aut(A)$ is a \hm\, then
$A\rtimes_{\Lambda}\Z^k$ can be embedded into an AF-algebra.

Suppose that there is a unital monomorphism $h:
C(X)\rtimes_{\Lambda}\Z^k\to A$ for some unital simple AF-algebra.
Suppose that $\tau$ is a tracial state on $A.$ Then $\tau\circ h$
gives a strictly positive $\Lambda$-invariant probability Borel
measure (every non-empty open set has a positive measure). In
other words, $\tau\circ h$ is a faithful $\Lambda$-invariant
tracial state. The main result of this paper is to show the
following theorem: Let $C$ be a unital AH-algebra, let $G$ be  a
finitely generated abelian group and let $\Lambda: G\to Aut(C)$ be
a \hm. Then the crossed product \CA\, $C\rtimes_{\Lambda} G$ can
be embedded into a unital simple AF-algebra if and only if $C$ has
a faithful $\Lambda$-invariant tracial state.  In particular, when
$A=C(X)$ and $G=\Z^k,$ this result implies that
$C(X)\rtimes_{\Lambda}\Z^k$ can be embedded into a unital simple
AF-algebra if and only if $X$ admits a strictly positive
$\Lambda$-invariant Borel probability  measure.

In \cite{Bn1}, \cite{M1}, as well as in \cite{Bn2}, a version of
so-called Berg's technique combined together with a version of
Rokhlin property played an important role. Let $C$ be a unital
separable amenable \CA\, and let $\af\in Aut(C)$ be an
automorphism. Suppose that one can formulate a right condition for
two unital monomorphisms $\phi_1, \phi_2: C\to A,$ where $A$ is a
unital simple AF-algebra, to be asymptotically unitarily
equivalent, i.e., there is a continuous path of unitaries $\{u_t:
t\in [0,\infty)\}$ of $A$ such that
$$
\lim_{t\to\infty}{\rm ad}\, u_t\circ \phi_1(a)=\phi_2(a)
$$
for all $a\in C.$ Then the technique developed in the above
mentioned articles can be used to prove that $C\rtimes_{\af}\Z$
can be embedded into a unital simple AF-algebra (see 5.5   of
\cite{Lnemb2}). In \cite{Lnasy}, we present a necessary and a
sufficient $K$-theoretical condition for $\phi_1$ and $\phi_2$
being asymptotically unitarily equivalent in the case that $C$ is
a unital AH-algebra. From this, we obtain a necessary and
sufficient condition for $C\rtimes_{\af}\Z$ to be embedded into a
unital simple AF-algebra in the case that $C$ is a unital
AH-algebra mentioned above (see 10.5 of \cite{Lnasy}). However, it
is not clear how to formulate a higher dimensional version of
Berg's technique.

On the other hand, one may write
$C\rtimes_{\Lambda}\Z^2=(C\rtimes_{\af_1}\Z)\rtimes_{\af_2}\Z$ for
some automorphisms $\af_1$ and $\af_2$ on $C.$  To prove a true
asymptotic unitary equivalence theorem for \CA s with the form
$C\rtimes_{\af_1}\hspace{-0.05in}\Z$ as above seems to require
much more general theory for the classification of amenable \CA s.
However, we try to pass it via  a special path. Given two
monomorphisms $\phi_1, \phi_2: C\rtimes_{\af_1}\Z\to A,$ where $A$
is a unital simple AF-algebra with a unique tracial state and
$K_0(A)=\rho_A(K_0(A)),$ i.e., $K_0(A)$ has no infinitesimal
elements, we will give a (necessary and sufficient) condition for
their induced maps $\phi_1^{(1)}, \phi_2^{(1)}$ from
$C\rtimes_{\af_1}\hspace{-0.05in}\Z$ to $B,$ a larger simple
AF-algebra, to be asymptotically unitarily equivalent. This will
lead to a monomorphism from $C\rtimes_{\Lambda}\Z^2$ into another
unital simple AF-algebra.

To establish such a result, as in the case of AH-algebra (see
\cite{Lnasy}), it requires some version of so-called Basic
Homotopy Lemma for the crossed products. It in turn requires an
approximate unitary equivalence theorem for the crossed products.
Our strategy is first to prove such a result. This requires a
non-commutative version of Berg's technique and a version of
Rokhlin property as well as a uniqueness theorem  for
monomorphisms from an AH-algebra. We also need an existence type
theorem to provide the needed  bott map. Combining these, we are
able to establish the required  homotopy lemma. Finally, to obtain
similar results  beyond $\Z^2$ action, one uses the induction.
Several arguments will be repeated.

The paper is organized as follows: Section 2 provides a number of
conventions, facts and reviews a couple of known results which
will be used several times in the proof. Section 3 is devoted to
the approximate unitary equivalence for monomorphisms from the
crossed products.  In section 4, we show that the prescribed bott
map can be constructed using some existing knowledge.  Section 5
establishes the required homotopy lemma. In section 6, we proved
the needed asymptotic unitary equivalence theorem for the crossed
products. In section 7, we present the theorem that
$C\rtimes_{\Lambda}\Z^2$ can be embedded into a unital simple
AF-algebra if and only if $C$ admits a faithful
$\Lambda$-invariant tracial state. Section 8 provides some
absorption lemma. Finally, in section 9, we use the induction to
prove the main embedding theorem.

\section{Preliminaries}

\begin{NN}

{\rm  Let $A$ be a \CA. Denote by $T(A)$ the tracial state space of $A$ and
denote by $Aff(T(A))$ the space of all real affine continuous functions on $T(A).$
Denote by $\rho_A: K_0(A)\to Aff(T(A))$ the positive \hm\, defined by $\rho_A([p])=\tau\otimes Tr(p)$ for all $\tau\in T(A)$ and projections $p\in M_k(A),$ where $Tr$ is the (non-normalized) standard trace.

}

\end{NN}

\begin{NN}

{\rm
Let $U(A)$ be a unitary group of $A.$ Denote by $Aut(A)$ the automorphism group.
If $u\in U(A),$ denote by ${\rm ad}\,u$ the inner automorphism defined by
${\rm ad}\, u(a)=u^*au$ for all $a\in A.$

}

\end{NN}

\begin{NN}

{\rm A \CA\, $A$ is an AH-algebra if $A=\lim_{n\to\infty}(A_n, \psi_n),$ where
each $A_n$ has the form $P_nM_{k(n)}(C(X_n))P_n,$ where $X_n$ is a finite CW complex (not necessarily
connected) and $P_n\in M_{k(n)}(C(X_n))$ is a projection.

We use $\psi_{n, \infty}: A_n\to A$ for the induced \hm.
Note that every separable commutative \CA\, is an AH-algebra. AF-algebras and  $A\T$-algebras
are AH-algebras.

}
\end{NN}

\begin{NN}

{\rm
Denote by ${\cal N}$ the class of separable amenable \CA s which satisfies the Universal Coefficient Theorem
(\cite{RS}).

}

\end{NN}

\begin{NN}

{\rm Let $X$ be a compact metric space. Denote by $Homeo(X)$ the group of all homeomorphisms on $X.$ }

\end{NN}
\begin{NN}

{\rm
Let $A$ and $B$ be two \CA s and let $L_1, L_2: A\to B$  be two maps. Let $\ep>0$ and
${\cal F}\subset A$ be a subset.
We write
$$
L_1\approx_{\ep} L_2\,\,\,\text{on}\,\,\, {\cal F},
$$
if
$$
\|L_1(a)-L_2(a)\|<\ep\rforal a\in {\cal F}.
$$

We say map $L_1$ is $\ep$-${\cal F}$-multiplicative if
$$
\|L_1(ab)-L_1(a)L_1(b)\|<\ep\rforal a, b\in {\cal F}.
$$

}
\end{NN}

\begin{df}\label{torus}
{\rm
Let $A$ and $C$ be two unital \CA s and let $\phi_1, \phi_2: C\to A$ be two unital monomorphisms.
Define the mapping torus
$$
M_{\phi_1,\phi_2}=\{f\in C([0,1], A): f(0)=\phi_1(c),\,\,\, f(1)=\phi_2(c)\,\,\,{\rm for\,\,\,some }\,\,\, c\in C\}.
$$
One
obtains an exact sequence:
\beq\label{Dmt-2}
0\to SA \stackrel{\imath}{\to} M_{\phi, \psi}\stackrel{\pi_0}{\to}
C\to 0.
\eneq
Denote by $\pi_t: M_{\phi_1, \phi_2}\to C$ the point-evaluation at $t$ ($t\in [0,1]$).

Suppose that $C$ is a separable amenable \CA.  From (\ref{Dmt-2}),
one obtains an element in $Ext(C,SA).$ In this case we identify
$Ext(C,SA)$ with $KK^1(C,SA)$ and $KK(C,A).$

Suppose that  $[\phi]=[\psi]$ in $KK(C,A)$ and $C$ satisfies the
Universal Coefficient Theorem, using Dadarlat-Loring's notation,
one has the following splitting exact sequence:
\beq\label{eta-1-}
0\to \underline{K}(SA)\,{\stackrel{[\imath]}{\to}}\,
\underline{K}(M_{\phi_1,\phi_2})\,{\stackrel{[\pi_0]}{\rightleftarrows}}_{\theta}
\,\,\underline{K}(C)\to 0.
\eneq
In other words there is  $\theta\in Hom_{\Lambda}(\underline{K}(C), \underline{K}(M_{\phi_1,\phi_2}))$ such that
$[\pi_0]\circ \theta=[\rm id_C].$ In particular, one has a
monomorphism $\theta|_{K_1(C)}: K_1(A)\to K_1(M_{\phi, \psi})$ such that
$[\pi_0]\circ \theta|_{K_1(A)}=({\rm id}_C)_{*1}.$
Thus, one may write
\beq\label{eta-1}
K_1(M_{\phi_1, \phi_2})=K_0(A)\oplus K_1(C).
\eneq

 }

\end{df}

\begin{df}\label{drot}

{\rm Suppose that $T(A)\not=\emptyset.$  Let $u\in M_l(M_{\phi,
\psi})$ be a unitary which is a piecewise smooth function on
$[0,1].$ For each $\tau\in T(A),$ we denote by $\tau$ for the
trace $\tau\otimes Tr$ on $M_l(A),$ where $Tr$ is the standard
trace on $M_l. $ Define
\beq\label{Dr-1}
R_{\phi,\psi}(u)(\tau)={1\over{2\pi i}}\int_0^1
\tau({du(t)\over{dt}}u(t)^*)dt.
\eneq
$R_{\phi, \psi}(u)(\tau)$ is real for every $\tau\in T(A).$

Suppose that $[\phi_1]=[\phi_2]$ in $KL(C,A).$ We also assume
that
\beq\label{Dr-2}
\tau(\phi_1(c))=\tau(\phi_2(c))\rforal c\in C\andeqn \tau\in T(A).
\eneq

Exactly as in section 2 of \cite{KK2}, one has the following
statement:}

{\it When $[\phi_1]=[\phi_2]$ in $KL(C,A)$ and {\rm (\ref{Dr-2})}
holds, there exists a \hm\,
$$
R_{\phi, \psi}: K_1(M_{\phi, \psi})\to Aff(T(A))
$$
defined by
$$
R_{\phi, \psi}([u])(\tau)={1\over{2\pi i}}\int_0^1
\tau({du(t)\over{dt}}u(t)^*)dt.
$$
}

 {\rm We will call $R_{\phi, \psi}$ the {\it rotation map} for the pair
$\phi$ and $\psi.$ }

{\rm Suppose also that $\tau\circ \phi_1=\tau\circ \phi_2$ for all $\tau\in T(A).$  Then  one obtains the \hm\,
\beq\label{eta-2}
R_{\phi,\psi}\circ \theta|_{K_1(C)}: K_1(C)\to Aff(T(A)).
\eneq

To keep the same notation as in \cite{KK2}, we write $${\tilde
\eta}_{\phi_1,\phi_2}=0,$$ if $R_{\phi,\psi}\circ \theta=0,$ i.e.,
$\theta(K_1(C))\in {\rm ker}R_{\phi_1,\phi_2}$ for some such
$\theta.$ Thus, $\theta$ also gives the following:
$$
{\rm ker}R_{\phi,\psi}={\rm ker}\rho_A\oplus K_1(C).
$$

Under the assumption that $[\phi_1]=[\phi_2]$ in $KK(C,A)$ and $\tau\circ \phi_1=\tau\circ \phi_2 $
for all $\tau\in T(A),$ and if, in addition, $K_i(C)$ is torsion free, such $\theta$ exists
whenever $\rho_A(K_0(A))=R_{\phi_1,\phi_2}(K_1(M_{\phi_1,\phi_2}))$ and
the following splits:
$$
0\to {\rm ker}\rho_A\to {\rm ker}R_{\phi_1, \phi_2}\to K_1(A)\to 0.
$$

For further information about rotation maps, see section 2 of \cite{KK2} and section 3 of \cite{Lnasy}.

}

\end{df}

\begin{NN}\label{ddbot}
{\rm Let $A$ be a unital amenable \CA, let $B$ be a unital \CA\,
and let $h: A\to B$ be a \hm. Suppose that $v\in U(B).$  We will
refer to \cite{Lnhomp} for the definition of $\text{Bott}(h, v)$ and $\text{bott}_1(h, v).$

Given a finite subset ${\cal P}\subset \underline{K}(A),$ there
exists a finite subset ${\cal F}\subset A$ and $\dt_0>0$ such that
$$
\text{Bott}(h, v)|_{\cal P}
$$
is well defined, if
$$
\|[h(a),\, v]\|=\|h(a)v-vh(a)\|<\dt_0\tforal a\in {\cal F}
$$
(see 2.10 of \cite{Lnhomp}).

There is $\dt_1>0$ (\cite{Lo}) such that
$\text{bott}_1(u,v)$ is well defined for any pair of unitaries $u$ and $v$ such that
$\|[u,\, v]\|<\dt_1.$ As in 2.2 of \cite{ER}, if $v_1,v_2,...,v_n$ are unitaries such that
$$
\|[u, \, v_j]\|<\dt_1/n,\,\,\,j=1,2,...,n,
$$
then
$$
\text{bott}_1(u,\,v_1v_2\cdots v_n)=\sum_{j=1}^n\text{bott}_1 (u,\, v_j).
$$


We will also use  ${\boldsymbol{\bt}}$ for the usual \hm\, from
$K_1(A)$ to $K_0(A\otimes C(\T)).$

See section 2 of \cite{Lnhomp} for the further information.
}

\end{NN}

\begin{df}\label{U0}
{\rm Denote by ${\cal U}$ throughout this paper the universal
UHF-algebra ${\cal U}=\otimes_{n\ge1}M_n.$

Let $\{e_{i,j}^{(n)}\}$ be the canonical matrix units for $M_n.$
Let $u_n\in M_n$ be the unitary matrix such that ${\rm ad}\,
u_n(e_{i,i}^{(n)})=e_{i+1,i+1}^{(n)}$ (modulo $n$). Let
$\sigma=\otimes_{n\ge 1} {\rm ad}\, u_n \in Aut({\cal U})$ be the
shift (see for example Example 2.2 of \cite{Bn1}). A fact that we
will use in this paper is the following  {\it cyclic Rokhlin property} that $\sigma$
has: For any integer $k>0,$ any $\ep>0$ and any finite subset
${\cal F}\subset {\cal U},$  there exist mutually orthogonal
projections $e_1,e_2,...,e_k\in {\cal U}$ such that $
\sum_{i=1}^ke_i=1_{\cal U},$ $ \|xe_i-e_ix\|<\ep$ for all $ x\in
{\cal F}$ and $\sigma(e_i)=e_{i+1}, i=1,2,..., k\,\,\, {\rm
(}e_{k+1}=e_1{\rm )}.$

Denote ${\cal U}^k=\overbrace{{\cal  U}\otimes {\cal  U} \otimes \cdots
\otimes {\cal U} }^k.$ We note that ${\cal U}^k\cong {\cal U}.$

Throughout this paper, ${\cal U}$ and $\sigma$ will be as the above.
}
\end{df}

\begin{df}\label{U1}

{\rm It follows from \cite{M1} (and its proof) that
 there is a unital monomorphism $j: {\cal U}
\times_\sigma\Z\to {\cal U} .$  However, in this case, more is
true. First $[\sigma]=[{\rm id}|_{U}]$ in $KK({\cal U}, {\cal U})$
and $\tau=\tau\circ \sigma.$ Since $K_1( {\cal U})=\{0\},$
$K_1(M_{{\rm id}_{\cal U}, \sigma})=K_0({\cal U}).$ In
particular (see  \cite{KK2}), there exists a continuous path of
unitaries $\{v(t): t\in [0, \infty)\}$ of ${\cal U}$ such that
\beq\label{dU1}
\lim_{t\to\infty}v(t)^*av(t)=\sigma(a)\rforal a\in  {\cal U} .
\eneq
Therefore, as in \cite{M1}, there is a unital embedding
$\phi:{\cal U}\rtimes_\sigma\Z\to {\cal U}$ such that
\beq\label{dU2}
\tau\circ \phi=\tau.
\eneq

Define $\psi: {\cal U}\rtimes_{\sigma}\Z\to {\cal U}\otimes {\cal U}$ by
$\psi(a)=\phi(a)\otimes 1_{\cal U}$ for all $a\in {\cal U}$ and $\psi(u_\sigma)=\phi(u_\sigma)\otimes \phi(u_{\sigma}^*).$
Then $\psi$ is a unital monomorphism.
Denote by $s: {\cal U}\otimes {\cal U}\to {\cal U}$ an isomorphism with $s_{*0}={\rm id}_{K_0({\cal U}\otimes {\cal U})}.$
We define $\imath: {\cal U}\rtimes_{\sigma}\Z\to {\cal U}$ by $s\circ \psi.$

{\it In what follows, $\imath$  will be used without further explanation.}

}

\end{df}

\begin{lem}\label{Urok}
Let $C$ be a unital separable \CA\, and let $\af\in Aut(C)$ be an automorphism.
Let $B$ be a unital separable amenable \CA.
Suppose that $\phi: C\rtimes_\af\Z\to B$ is  a unital \hm\, such that
$\phi|_C$ is injective. Define a unital \hm\, $\phi^{(1)}: C\rtimes_\af\Z\to B\otimes {\cal U}$ defined
by $\phi^{(1)}(c)=\phi(c)\otimes 1_{\cal U}$ for all $c\in C$ and $\phi^{(1)}(u_\af)=\phi(u_\af)\otimes \imath(u_\sigma).$

Then $\phi^{(1)}$ is a monomorphism.
\end{lem}

\begin{proof}
Let $C_1=(C\otimes{\cal U})\rtimes_{\af\otimes \sigma}\Z.$ Then $\af\otimes \sigma$ is known
to have the cyclic Rokhlin property (see \ref{U0}). We view ${\cal U}$ as a \SCA\, of ${\cal U}\rtimes_{\sigma}\Z.$
Define $\psi: C_1\to B\otimes {\cal U}$ by
$\psi(c\otimes b)=\phi(c)\otimes \imath(b)$ for all $c\in C$ and $b\in {\cal U}$ and $\psi(u_{\af\otimes \sigma})=\phi(u_\af)\otimes \imath(u_\sigma).$
Then $C_1|_{C\otimes {\cal U}}$ is injective. Since $\af\otimes \sigma$ has the cyclic Rokhlin property,
by Lemma 4.1 of \cite{Lnemb1}, $\psi$ is injective.
Let $j_0: C\rtimes_\af\Z\to C_1$ be defined by $j_0(c)=c\otimes 1_{\cal U}$ for $c\in C$ and $j_0(u_\af)=u_{\af\otimes \sigma}.$
Then $j_0$ is an embedding. Note $\phi^{(1)}=\psi\circ j_0.$

\end{proof}

\begin{thm}{\rm(Corollary 4.8 of \cite{Lncd})} \label{2QUni}
Let $C$ be a unital AH-algebra and let $\tau_0\in T(C)$ be a
faithful tracial state.

Let $\ep>0,$ ${\cal F}\subset C$ be a finite subset. There exists
$\dt>0,$ $\sigma>0,$ a finite subset ${\cal G}$ and a finite
subset ${\cal P}\subset \underline{K}(C)$ satisfying the
following:

For any unital simple \CA\, with tracial rank zero with a unique
tracial state $\tau\in T(A),$ any two unital $\dt$-${\cal
G}$-multiplicative \morp s $L_1, L_2: C\to A$ such that
\beq\label{2QU1}
[L_1]|_{\cal P}=[L_2]|_{\cal  P}\andeqn\\
\tau\circ L_i\approx_{\sigma} \tau_0\,\,\,\text{on}\,\,\,{\cal  G},
\eneq
there exists a unitary $w\in U(A)$ such that
\beq\label{2QU2}
{\rm ad}\, w\circ L_1\approx_{\ep} L_2\,\,\,\text{on}\,\,\,{\cal
F}.
\eneq

\end{thm}

\begin{proof}
It is clear that the general case can be reduced to the case that $C=PM_k(C(X))P,$
where $X$ is a finite-dimensional compact metric space and $P\in M_k(C(X))$ is a projection.
Since $\tau_0$ is fixed, the statement follows from Cor. 4.8 of \cite{Lncd}.

\end{proof}

The following is a direct application of a theorem Choi and Effros (\cite{CE}, see 5.10.10 of \cite{Lnbk}).

\begin{prop}\label{APP}
Let $A$ be a separable amenable \CA. For any $\ep>0$ and finite subset ${\cal F}\subset A,$ there is $\dt>0$ and
a finite subset ${\cal G}\subset A$ satisfying the following:

For any \CA\, $B$ and any self-adjoint linear map $L: A\to B$ which is $\dt$-${\cal G}$-multiplicative,
there exists an $\ep$-${\cal F}$-multiplicative \morp\, $\phi: A\to B$ such that
\beq\label{app1}
L\approx_{\ep} \phi\,\,\,\text{on}\,\,\, {\cal F}.
\eneq

If, furthermore, $A$ has a unital, one may assume that $\phi(1_A)$ is a projection.
\end{prop}

\section{Approximate unitary equivalence}

\begin{lem}\label{burgl1}

Let $C$ be a unital amenable separable \CA\,  and $\af\in Aut(C)$ be an automorphism.
For any $\ep>0,$ $L>0$ and finite subset ${\cal  F}\subset C,$ there is $\dt>0$ and  a finite subset ${\cal  G}\subset C$ satisfying the following:

Suppose that $A$ and $B$ are two unital amenable separable \CA s,
$\phi: C\rtimes_{\af}\Z\to A$ and $\psi: {\cal
U}\rtimes_{\sigma}\Z\to B$ are two unital monomorphisms. Fix
unitaries  $u_1\in A,$ $u_2\in B$ such that $u_1^*\phi(a)u_1=\phi(
\af(a))$ for all $a\in C$ and $u_2^*\psi(b)u_2=\psi(\sigma(b))$ for
all $b\in {\cal  U}.$ Suppose that there is a unitary $v\in U(A)$
and a continuous path of unitaries $\{v_t: t\in [0,1]\}\subset
A\otimes 1$ such that
\beq\label{bl-1}
&&v_0=1,\,\,\,v_1=\prod_{j=0}^{n-1}(u^{n-1-j})^*(v\otimes 1_B)u^{n-1-j}  ,\\
&&\|[v,\, \phi(a)]\|<\dt,\,\,\,\|[v_t,\, \phi(a)\otimes 1_B]\|<\dt\tforal a\in {\cal  G} \tand  t\in [0,1],\\
&&\tand Length(\{u_t\})\le L\andeqn (L+1)/n<\ep/2,
\eneq
where $u=u_1\otimes u_2,$
then there is a unitary $w\in A\otimes B$ such that
\beq\label{bl-2}
{\rm ad}\, w\circ\phi(a)\approx_{\ep} \phi(a)\tforal a\in {\cal
F}\andeqn w^*u(v\otimes 1_B)w\approx_{\ep} u,
\eneq
where we identify $a$ with $a\otimes 1$ for $a\in A.$
\end{lem}

\begin{proof}
Fix $\ep>0.$ We may assume that ${\cal  F}$ is in the unit ball of $C.$
Put $\dt={\ep\over{4n(n+1)}}.$
Let ${\cal  G}=\cup_{-2n\le j\le 2n} \af^j({\cal  F}).$ Suppose that $\{v_t\}$ satisfies the assumption.

There  are mutually orthogonal projections $q_1, q_2,...,q_n\in
\psi({\cal  U})$ with $\sum_{j=1}^n q_j=1_{\cal U}$ such that
\beq\label{bl1}
u_2^*q_ju_2=q_{j+1},\,\,\,j=1,2,....,n\,\,\,\hspace{0.3in} \text{(} q_{n+1}=q_1 \text{)}.
\eneq
Define $e_i=1\otimes q_i,$ $i=1,2,...,n+1$ with $e_{n+1}=e_1.$

In what follows, we will identify an element $a\in A$ with the element $a\otimes 1$ in $A\otimes B$
whenever it is convenient. So we may write that $v=v\otimes 1$ and note that
$$
v_1=((u^{n-1})^*vu^{n-1})((u^{n-2})^*vu^{n-2})((u^{n-3})^*v(u^{n-3})\cdots u^*vuv.
$$

There are $y_1, y_2,...,y_n\in \{v_t\}\subset A\otimes 1$ such that
\beq\label{bl2}
y_1=1,\,\,\, y_n=v_1\andeqn \|y_j-y_{j+1}\|<{L\over{n}}, \,\,\, j=1,2,...,n.
\eneq
Put $x_1=y_1^*, $ $x_2=y_1y_2^*,$ $x_3=y_2y_3^*,$..., $x_{n-1}=y_{n-2}y_{n-1}^*$ and $x_n=y_{n-1}y_n^*.$
Note that
\beq\label{bl3}
\|[\phi(a),\, x_j]\|<2\dt\rforal a\in {\cal  G},\,\,\,j=1,2,...,n.
\eneq
Put
$$
s_j=v^*u^*vu\cdots (u^{j-1})^*v^*u^{j-1},\,\,\,j=1,2,...,n-1
$$
Define
\beq\label{bl4}
\eta_1&=& u^{n-1}s_{n-1}^*x_1s_{n-1}(u^{n-1})^*\\
\eta_2&=&u^{n-2}s_{n-2}^*x_2s_{n-2}(u^{n-2})^*,\\
&\cdots&,\\
\eta_{n-2}&=& u^2 s_2^*x_{n-2}s_2(u^2)^*=u^2(u^*vu)vx_{n-2}v^*(u^*v^*u)(u^2)^*\\
\eta_{n-1}&=&uvx_{n-1}v^*u^*\\
\eta_n&=&x_n.
\eneq
Note that
\beq\label{bl5}
\phi(a)uvu^*&=&uu^*\phi(a)uvu^*=u\phi\circ \af(a)vu^*\\
&\approx_{\dt}& uv\phi\circ \af(a)u^*\\
&=&uvu^*u\phi\circ \af(a)u^*=uvu^*\phi(a)
\eneq
for all $a\in \cup_{-2(n-1)\le j\le 2(n-1)}\af({\cal  F}).$
It follows that
\beq\label{bl6}
\|[\phi(a), \, \eta_j]\|<2(n-j+1)\dt\rforal a\in \cup_{-2(n-j)\le j\le 2(n-j)}\af^j({\cal  F})
\eneq
and for $ j=1,2,...,n.$
Put
\beq\label{bl7}
&&d_1=e_1v\eta_1e_1,\,\,\, d_2=e_2v\eta_2e_2, \,\,\,d_3=u^*d_2u v\eta_3e_3,...\\
&&d_{j+1}=u^*d_{j}u v\eta_{j+1}e_{j+1},..., d_n=u^*d_{n-1}uv\eta_ne_n.
\eneq

Note that
\beq\label{bl8}
d_n&=&(u^2)^*d_{n-2}u^2 u^*v(uvx_{n-1}v^*u^*)uv\eta_ne_n\\
&=& (u^2)^*d_{n-2}u^2 u^*vuvx_{n-1}x_ne_n\\
&=&(u^3)^*d_{n-3}u^3  (u^2)^*v
(u^2 u^*vuvx_{n-2}v^*(u^*vu)(u^2)^*)(u^2u^*vuvx_{n-1}x_ne_n)\\
&=&(u^3)^*d_{n-3}u^3  (u^2)^*vu^2 u^*vuvx_{n-2}x_{n-1}x_ne_n\\
&\cdots&\\
&=&v_1x_1x_2\cdots x_ne_n=e_n
\eneq

Define
\beq\label{bl9}
Z=\sum_{j=1}^n  d_j.
\eneq
Note that $Z\in U(A\otimes B).$
We estimate that
\beq\label{bl10}
\phi(a)d_1&=&e_1\phi(a)ve_1\eta_1e_1\\
&\approx_{\dt}& e_1v\phi(a)\eta_1e_1\\
&\approx_{2n\dt}& e_1v\eta_1e_1\phi(a)=d_1\phi(a)
\eneq
for all $a\in {\cal  G}.$
Suppose that we have shown
\beq\label{bl11}
\phi(a)d_j\approx_{(j(2n+3)-j(j+1))\dt} d_j\phi(a)\rforal a\in \cup_{-2(n-j)\le i\le 2(n-j)}\af^i({\cal  F}).
\eneq
\beq\label{bl11+}
\phi(a) d_{j+1}&=&\phi(a) u^*d_juv\eta_je_{j+1}\\
&=&u^*\phi\circ \af^{-1}(a)d_juv\eta_je_{j+1}\\
&\approx_{(j(2n+3)-j(j+1))\dt}& u^*d_j\phi\circ \af^{-1}(a)uv_1\eta_je_{j+1}\\
&=&u^*d_ju\phi(a)v\eta_je_{j+1}\\
&\approx_{\dt} &u^*d_juv\phi(a)\eta_je_{j+1}\\
&\approx_{2(n-j)\dt}& u^*d_juv\eta_je_{j+1}\phi(a)
\eneq
for all $a\in \cup_{-2(n-j-1)\le i\le 2(n-j-1)}\af^i({\cal  F}).$
In other words,
\beq\label{bl12}
\phi(a)d_{j+1}\approx_{(j+1)(2n+3)-(j+1)(j+2)\dt} d_{j+1}\phi(a)
\eneq
for all $a\in  \cup_{-2(n-j-1)\le i\le 2(n-j-1)}\af^i({\cal  F}).$
It follows that
\beq\label{bl13}
\|[Z, \, \phi(a)]\|<2n(n+1)\dt \rforal a\in {\cal  F}.
\eneq
Note that
\beq\label{bl13-}
\|\eta_j-1\|<{L/n}<\ep/2,\,\,\,j=1,2,...,n.
\eneq

We then compute that
\beq\label{bl14}
Z^*uZ&=&d_n^*e_nue_1d_1+ \sum_{j=1}^{n-1} d_{j}^*e_{j}ue_{j+1}d_{j+1}\\
&=& e_nuv\eta_1e_1+\sum_{j=1}^{n-1} u (u^*d_j^*ue_{j+1}d_{j+1})\\
&=& e_nuv\eta_1e_1+\sum_{j=1}^{n-1} u(u^*d_j^*u e_{j+1}u^*d_{j}uv\eta_je_{j+1})\\
&=&e_nuv\eta_1e_1+\sum_{j=1}^{n-1}ue_{j+1}v\eta_je_{j+1}\\\label{bl14+}
&\approx_{\ep}& e_nuve_1+\sum_{j=1}^{n-1}uve_{j+1}=uv\,\,\,\,\,\,\,\,\,\,\,\,\,\,\,\,\,\,{\rm (applying (\ref{bl13-}))}.
\eneq
Put $w=Z^*.$ It follows from (\ref{bl13}) that
\beq\label{bl15}
{\rm ad}\, w\circ \phi(a)\approx_{\ep}\phi(a)\tforal a\in {\cal  F}.
\eneq
Moreover, by (\ref{bl14}),
\beq\label{bl16}
w^*(u_1v\otimes u_2)w=Z(u_1v\otimes u_2)Z^*\approx_{\ep} u.
\eneq

\end{proof}

\begin{lem}\label{2burg}
Let $C$ be a unital AH-algebra and let $\af\in Aut(C)$ be an automorphism.
Suppose that $A$  and $B$ are two unital simple \CA s with tracial
rank zero, $\phi: C\rtimes_{\af}\Z\to A$ and $\psi: {\cal
U}\rtimes_{\sigma}\Z\to B$ are two unital monomorphisms.

For any $\ep>0$ and any finite subset ${\cal  F}\subset C,$ there is $\dt>0,$ a
finite subset
${\cal  G}\subset C$ and a finite subset ${\cal  P}\subset \underline{K}(C)$
satisfying the following:

 Fix
unitaries  $u_1\in A,$ $u_2\in B$ such that $u_1^*\phi(a)u_1=\phi(
\af(a))$ for all $a\in C$ and $u_2^*\psi(b)u_2=\psi(\sigma(b))$
for all $b\in {\cal  U}.$ Suppose that there is a unitary $v\in
U(A)$ such that
\beq\label{2bl-1}
\|[v,\, \phi(a)]\|<\dt,\,\,\,\tforal a\in {\cal  G} \tand
 {\rm Bott}(\phi, \, v)|_{\cal  P}=0,
\eneq
then there is a unitary $w\in A\otimes B$ such that
\beq\label{2bl-2}
w^*(\phi(a)\otimes 1_B)w\approx_{\ep} \phi(a)\otimes 1_B \tforal
a\in {\cal  F}\andeqn w^*u(v\otimes 1_B)w\approx_{\ep} u,
\eneq
where $u=u_1\otimes u_2.$
\end{lem}

\begin{proof}
Fix a finite subset ${\cal  F}\subset C$ and $\ep>0.$ Define
$\phi': C\to A\otimes B$ by $\phi'(c)=\phi(c)\otimes 1_B$ for all
$c\in C.$

 Let $\dt_1>0$ (in place of $\dt$)  and ${\cal
G}_1\subset C$ (in place of ${\cal  G}$) be a finite subset
required by \ref{burgl1} for $L=2\pi+1.$

Let $n\ge 1$ be an integer such that
${2\pi+1\over{n}}<\min\{\ep/2, \dt_1/2\}.$

Let $\dt_2>0$ (in place of $\dt$),  ${\cal  G}_2\subset C$ (in
place of ${\cal  G}$) be a finite subset and let ${\cal  P}_1\subset
\underline{K}(C)$ be a finite subset required by 17.9 of \cite{Lnhomp}
corresponding to $\dt_1/2$ (in place of $\ep$), ${\cal  G}_1$ (in
place of ${\cal  F}$) and $\phi'|_C.$

Without loss of generality, to simplify notation, we may assume
that ${\cal  F},$ ${\cal  G}_1$ and ${\cal  G}_2$ are all in the unit
ball of $C.$  We may further assume that $\dt_2<\min\{\dt_1/2,
\ep/2\}.$

Define
\beq\label{2bl-3}
v_1= \prod_{j=0}^{n-1}(u^{n-1-j})^*(v\otimes 1)u^{n-1-j},
\eneq
where $u=u_1\otimes u_2.$ Note that $v_1\in A\otimes 1.$
 Define
$$
{\cal  G}=\cup_{k=-2n}^{2n}\af^k({\cal  G}_2)\andeqn {\cal
P}=\cup_{k=-2n}^{2n}[\af^k]({\cal  P}_1).
$$

Let $\dt>0$ so that $\dt<\dt_2/2n^2.$

Suppose that
\beq\label{2bl-4}
\|[\phi(c),\, v]\|<\dt\tforal c\in {\cal  G}\andeqn
\text{Bott}(\phi,\, v)|_{\cal  P}=0.
\eneq
We may also assume that $\text{Bott}(h', v')|_{\cal  P}$ is well
defined provided that
\beq\label{2bl-4+}
\|[h',\, v']\|<n^2 \dt<\dt_2/2
\eneq
for any unital \hm\, $h'$ from $C$ and any unitary $v'.$

Note we have
\beq\label{2bl-4++}
\text{Bott}(\phi', \, v\otimes 1_B)|_{\cal  P}=0.
\eneq

Note that,
\beq\label{2bl-5}
\phi'(a)u^*(v\otimes 1_B) u &=& u^*\phi'(\af^{-1}(a))(v\otimes 1_B)u\\
&\approx_{\dt}& u^*(v\otimes 1_B)\phi'(\af^{-1}(a))u\\
&=&u^*(v\otimes 1_B)u\phi'(a)
\eneq
for all $a\in \cup_{k=-2n+1}^{2n-1}\af^k({\cal  G}_2).$ It follows
that
\beq\label{2bl-6}
\|[\phi'(a), \, v_1]\|<n\dt<\dt_2/2n\rforal a\in {\cal  G}_2.
\eneq

By (\ref{2bl-4++}), we have that
\beq\label{2bl-7}
\text{Bott}(\phi',\, (u^j)^*(v\otimes 1_B)u^j)|_{{\cal  P}_1}
=\text{Bott}(\phi', v\otimes 1_B)|_{[\af^{-j}]({\cal  P}_1)}=0,\\
\eneq
$j=0,1,2,...,n-1.$ It follows that
\beq\label{2bl-8}
\text{Bott}(\phi',\,
v_1)|_{{\cal P}_1}=\sum_{j=1}^{n-1}\text{Bott}(\phi',\,(u^j)^*(v\otimes
1_B)u^j)|_{{\cal  P}_1}=0
\eneq
The above is computed in $A\otimes 1.$  It follows from the Basic
Homotopy Lemma 17.9 of \cite{Lnhomp} that there exists a continuous
path of unitaries $\{v(t): t\in [0,1]\}\subset U(A\otimes 1)$ such
that
\beq\label{2bl-9}
v(0)=1,\,\,\, v(1)=v_1\andeqn \|[\phi'(c),\,
v(t)]\|<\dt_1/2\rforal a\in {\cal  G}_1.
\eneq
Moreover
\beq\label{2bl-10}
\text{Length}(\{v(t)\})\le 2\pi+1=L.
\eneq

Thus \ref{burgl1} applies.

\end{proof}

\begin{cor}\label{4Tuni}
Let $C$ be a unital AH-algebra, let $\af\in Aut(C)$ be an
automorphism and let $A\cong A\otimes {\cal U}$ be a unital simple AF-algebra with a
unique tracial state $\tau$ and $K_0(A)=\rho_A(K_0(A)).$
 Suppose that $\phi_1, \phi_2: C\rtimes_{\af}\Z\to A$
are two unital monomorphisms.
For any $\ep>0$ and any finite subset ${\cal  F}\subset
C\rtimes_\af\Z,$ there exists $\dt>0,$ a finite subset ${\cal
G}\subset C$ and a finite subset ${\cal  P}\subset K_1(C)$
satisfying the following:

Suppose that there exists a unitary $z\in A$ such that
\beq\label{4Tu-1+}
{\rm ad}\, z\circ \phi_1(a)\approx_{\dt}\phi_2(a)\tforal a\in
{\cal  G}
\eneq
and, for any $x\in {\cal  P},$
\beq\label{4Tu-1++}
{\rm{bott}}_1(\phi_2,\, \phi_2(u_\af)^*z^*\phi_1(u_\af)z)(x)=0.
\eneq
Then there exists a unitary $w \in U(A\otimes {\cal U})$ such
that
\beq\label{4Tu-2}
{\rm ad}\, w\circ
\phi_1^{(1)}\approx_{\ep}\phi_2^{(1)}\,\,\,\text{on}\,\,\,{\cal  F},
\eneq
where $\phi_i^{(1)}: C\rtimes_\af\Z\to A\otimes{\cal U}$ is defined by
$\phi_i^{(1)}(c)=\phi_i(c)\otimes 1_{ {\cal U}}$ and
$\phi_i^{(1)}(u_\af)=\phi_i(u_\af)\otimes \imath(u_\sigma),$ $i=1,2.$

\end{cor}

\begin{proof}

Let $\ep>0$ and ${\cal  F}\subset C\rtimes_\af\Z$ be a finite subset.
Without loss of generality, we may assume that ${\cal  F}={\cal
F}_0\cup\{u_\af\},$ where ${\cal  F}_0\subset C$ is a finite
subset.

 Let $\dt>0,$ ${\cal  G}\subset C$ be a finite subset and ${\cal
P}_1\subset \underline{K}(C)$ be a finite subset required by
\ref{2burg} for $\ep/2$ and ${\cal  F}_0$ above (with $\phi$ being replaced by $\phi_2,$ $B={\cal U}$ and $\psi$ being replaced by $\imath$). Without loss of
generality, we may assume that ${\cal  G}$ is in the unit ball of
$C$ and $\dt<\ep/2.$

Since $K_0(A)$ is divisible and torsion free, there is a finite
subset ${\cal  P}\subset K_1(C)$ such that
\beq\label{4tu-1}
\text{Bott}(\phi_2,\, v)|_{{\cal  P}_1}=0
\eneq
provided that
\beq\label{4tu-2}
\text{bott}_1(\phi_2, v)|_{\cal  P}=0
\eneq
and $v$ is a unitary and $\text{Bott}(\phi_2,\, v)|_{{\cal  P}_1}$
is well-defined.

Now suppose that (\ref{4Tu-1+}) and (\ref{4Tu-1++}) hold for above
$\dt,$ ${\cal  G}$ and ${\cal  P}.$

Define
\beq\label{4tu-3}
v=\phi_2(u_\af)^*z^*\phi_1(u_\af)z.
\eneq
Then
\beq\label{4tu-4}
\|[\phi_2(c),\, v]\|<\dt\rforal c\in {\cal  G} \andeqn
\text{Bott}(\phi_2, \, v)|_{{\cal  P}_1}=0.
\eneq

It follows from \ref{2burg} that there exists a unitary $w_1\in
A\otimes {\cal U}$ such that
\beq\label{4tu-5}
{\rm ad}\, w_1\circ \phi_2^{(1)}\approx_{\ep/2}
\phi_2^{(1)}\,\,\,\text{on}\,\,\, {\cal  F}_0 \andeqn\\
w_1^*\phi_2^{(1)}(u_\af)(v\otimes 1_{
{\cal U}})w_1\approx_{\ep/2}\phi_2^{(1)}(u_\af).
\eneq
Note that
\beq\label{4tu-6}
\phi_2^{(1)}(u_\af)(v\otimes 1_{  {\cal U}})&=& \phi_2(u_\af)v\otimes
\imath(u_\sigma)\\
&\approx_{\dt}&z^*\phi_1(u_\af)z\otimes \imath(u_\sigma).
\eneq
Put $w=(z\otimes 1_{\cal U})w_1.$  Then, by (\ref{4tu-5}),
(\ref{4Tu-1+}) and (\ref{4tu-6}),
\beq\label{4tu-7}
{\rm ad}\, w\circ \phi_1^{(1)}\approx_{\ep}
\phi_2^{(1)}\,\,\,\text{on}\,\,\,{\cal  F}_0\andeqn {\rm ad}\,w\circ
\phi_1^{(1)}(u_\af)\approx_{\ep}\phi_2^{(1)}(u_\af).
\eneq

\end{proof}

\begin{cor}\label{3Tuni}
Let $C$ be a unital AH-algebra, let $\af\in Aut(C)$ be an
automorphism and let $A\cong A\otimes {\cal U}$ be a unital simple AF-algebra with a
unique tracial state $\tau$ and $K_0(A)=\rho_A(K_0(A)).$  Suppose that $\phi_1, \phi_2: C\rtimes_{\af}\Z\to A$
are two unital monomorphisms
and suppose that there exists a sequence of unitaries $z_n\in A$ such
that
\beq\label{3Tu-1+}
\lim_{n\to\infty}{\rm ad}\, z_n\circ \phi_1(a)=\psi_2(a)\tforal
a\in C
\eneq
and, for any $x\in K_1(C)$ and
\beq\label{3Tu-1++}
{\rm{bott}}_1(\phi_2,\, \phi_2(u_\af)^*z_n^*\phi_1(u_\af)z_n)(x)=0
\eneq
for all sufficiently large $n.$

  Then there exists a sequence of
unitaries $\{w_n\}\subset U(A\otimes   {\cal U})$ such that
\beq\label{3Tu-2}
\lim_{n\to\infty}{\rm ad}\, w_n\circ \phi_1^{(1)}(a)=\phi_2^{(1)}(a)\tforal
a\in C\rtimes_{\af}\Z,
\eneq
where $\phi_i^{(1)}: C\rtimes_\af\Z\to A\otimes {\cal U}$ is defined
by $\phi_i^{(1)}(u_\af)=\phi_i(u_\af)\otimes \imath(u_\sigma),$ $i=1,2.$
\end{cor}

\begin{thm}\label{Tuni}
Let $C$ be a unital AH-algebra with  $K_1(C)=\{0\},$ let $\af\in
Aut(C)$ be an automorphism and let $A\cong A\otimes {\cal U}$ be a unital simple
AF-algebra with a unique tracial state $\tau$ and
$K_0(A)=\rho_A(K_0(A)).$  Suppose that $\phi_1,
\phi_2: C\rtimes_{\af}\Z\to A$ are two unital monomorphisms such
that
\beq\label{Tu-1}
\tau\circ \phi_1=\tau\circ \phi_2.
\eneq

 Then there exists a sequence of unitaries
$\{w_n\}\subset U(A\otimes {\cal U})$ such that
\beq\label{Tu-2}
\lim_{n\to\infty}{\rm ad}\, w_n\circ \phi_1^{(1)}(a)=\phi_2^{(1)}(a)\tforal a\in C\rtimes_{\af}\Z.
\eneq
where $\phi_i^{(1)}: C\rtimes_{\af}\Z\to A\otimes {\cal U}$ by
$\phi_i^{(1)}(c)=\phi_i(c)\otimes 1_{\cal U}$ for all $c\in C$ and
$\phi_i^{(1)}(u_\af)=\phi_i(u_\af)\otimes \imath(u_\sigma),$ $i=1,2.$

\end{thm}

\begin{proof}
Since $K_0(A)=\rho_A(K_0(A)),$ (\ref{Tu-1}) implies that
\beq\label{Tu-3}
(\phi_1|_{C})_{*0} =(\phi_2|_{C})_{*0}.
\eneq
Since $K_1(A)=\{0\}$ and $K_0(A)$ is divisible, (\ref{Tu-3}) implies
that
\beq\label{Tu-4}
[\phi_1|_{C}]=[\phi_2|_{C}]\,\,\, \text{in}\,\,\, KK(C, A).
\eneq
Combining (\ref{Tu-4}) and (\ref{Tu-1}), by applying \ref{2QUni},
there exists a sequence of unitaries $\{z_n\}\in U(A)$ such that
\beq\label{Tu-5}
\lim_{n\to\infty}{\rm ad} z_n\circ \phi_1(c)=\phi_2(c)\tforal c\in C.
\eneq

Since $K_1(C)=0,$
\beq\label{Tu-6}
\text{bott}_1(\phi_2,\, \phi_2(u_\af)^*z_n^*\phi_1(u_\af)z_n)=0.
\eneq

We then apply \ref{3Tuni}.
\end{proof}

\begin{thm}\label{2Tuni}
Let $C$ be a unital AH-algebra, let $\af\in Aut(C)$ be an automorphism and let $A\cong A\otimes {\cal U}$ be a unital simple
AF-algebra with a unique tracial state $\tau$ and
$K_0(A)=\rho_A(K_0(A)).$  Suppose that $\phi_1,
\phi_2: C\rtimes_{\af}\Z\to A$ are two unital monomorphisms such
that
\beq\label{2Tu-1}
\tau\circ \phi_1=\tau\circ \phi_2.
\eneq
Suppose also that
\beq\label{2Tu-1+}
R_{\phi_1\circ j_0, \phi_2\circ j_0}(K_1(C))\subset
\rho_A(K_0(A)),
\eneq
where $j_0: C\to C\rtimes_\af \Z$ is the embedding.

  Then
there exists a sequence of unitaries $\{w_n\}\subset U(A\otimes
{\cal U})$ such that
\beq\label{2Tu-2}
\lim_{n\to\infty}{\rm ad}\, w_n\circ \phi_1^{(1)}(a)=\phi_2^{(1)}(a)\tforal a\in C\rtimes_{\af}\Z.
\eneq
where $\phi_i^{(1)}: C\rtimes_{\af}\Z\to A\otimes {\cal U}$ by
$\phi_i^{(1)}(c)=\phi_i(c)\otimes 1$ for all $c\in C$ and
$\phi_i^{(1)}(u_\af)=\phi_i(u_\af)\otimes \imath(u_\sigma),$ $i=1,2.$ Then
\end{thm}

\begin{proof}
Exactly as in the proof of \ref{Tuni}, we have
\beq\label{2Tu-3}
[\phi_1|_C]=[\phi_2|_C]\,\,\,\text{in}\,\,\,KK(C, A).
\eneq

Define
$$
M_{\phi_1, \phi_2}=\{f\in C([0,1], A): f(0)=\phi_1(a)\andeqn f(1)=\phi_2(a)\, \,\,\text{for\,\,\,some}\,\,\, a\in C\}.
$$
Then the condition (\ref{2Tu-1+}) implies that
\beq\label{2Tu-3+}
R_{\phi_1, \phi_2}(K_1(M_{\phi_1, \phi_2}))=\rho_A(K_0(A)).
\eneq
Then, by combining (\ref{Tu-4}), by applying Lemma 11.3 of \cite{Lnasy} and Theorem 10.7 of \cite{Lnasy}, one obtains a continuous path of unitaries
$\{z_t: t\in [0, \infty)\}\subset U(A)$ such that
\beq\label{2Tu-4}
\lim_{t\to\infty}{\rm ad}\, z_t\circ \phi_1(a)=\phi_2(a)\tforal a\in C.
\eneq
As in the proof of \ref{Tuni}, put
\beq\label{2Tu-5}
v_t=\phi_2(u_\af)z_t^*\phi_1(u_\af)z_t\tforal t\in [0, \infty).
\eneq
 Then one obtains
\beq\label{2Tu-6}
\lim_{t\to\infty}\|[\phi_2(a), \, v_t]\|=0\tforal a\in C.
\eneq
In particular, for a fixed finite subset ${\cal  P}\subset K_1(C),$
\beq\label{2Tu-7}
\text{bott}_1(\phi_2,\,  v_t)|_{\cal  P}=\text{bott}_1(\phi_2,\,
v_{t'})|_{\cal  P}
\eneq
for all sufficiently large $t, t'.$
By the Exel formula (3.6 of \cite{Lnasy}),  for each $z\in U(M_k(C)),$
one has
\beq\label{2Tu-8}
\lim_{t\to\infty}\text{bott}_1(\phi_2, v_t)([z])(\tau)=
\lim_{t\to\infty}\tau(\log(V_t^*\phi(z)V_t\phi(z)^*)=0,
\eneq
where $ V_t={\rm diag}(\overbrace{ v_t, v_t, ...,v_t}^k).$
Combining (\ref{2Tu-7}) and (\ref{2Tu-8}), one must have, for any
$x\in  {\cal  P},$
\beq\label{2Tu-9}
\text{bott}_1(\phi_2, v_1)(x)(\tau)=0.
\eneq
Since $K_0(A)=\rho_A(K_0(A)),$ this implies that
\beq\label{2Tu-10}
\text{bott}_1(\phi_2, v_t)(x)=0
\eneq
for all $x\in {\cal  P}$ for all large $t.$ Then \ref{3Tuni}
applies.

\end{proof}

\section{Bott maps}

\begin{thm}\label{Kemb1}
Let $C$ be a unital separable amenable \CA\, in ${\cal  N}$  which
has a unital embedding $j: C\to A$ for some unital simple
AF-algebra $A.$

Then there is a unital separable simple \CA\, $B\in {\cal  N}$
with tracial rank zero and a unital monomorphism $\phi: C\to B$
such that $\phi_{*1}=\gamma\circ \kappa,$ where $\kappa:
K_1(C)\to K_1(C)/\rm{Tor}(K_1(C))$ is the quotient map and
$\gamma: K_1(C)/{\rm Tor}(K_1(C))\to K_1(B)$ is injective.
Moreover, we may assume that $B$ has a unique tracial state,
if $A$ has a unique tracial state.

\end{thm}

\begin{proof}
Without loss of generality, we may assume that $K_1(C)$ is torsion
free, by replacing $C$ by  $C\otimes {\cal  U}.$

Let $\tau\in T(A)$ be a tracial state of $A.$ There is a unital
simple AF-algebra $A_0$ with a unique tracial state and
$K_0(A_0)=\rho_A(K_0(A_0))=\D,$ where $\D$ is a subring of $\R$
which is also a countable divisible dense subgroup of $\R$
containing $\tau(K_0(A))$ with $[1_{A_0}]=1.$ There is a unital
monomorphism $\psi_1: A\to A_0.$ Therefore, by replacing $j$ by
$\psi_1\circ j,$ we may assume that $A=A_0.$
We may also assume that $A\otimes {\cal U}=A.$

 Let $C_1=C \otimes A$ and
 $C_2=C_1\otimes A\otimes M_2,$...,$C_{n+1}=C_n\otimes A\otimes M_{n},$
$n=1,2,...,.$ We will identify $A$ with $A\otimes M_n$ whenever it
is convenient. In particular, $C_{n+1}\cong C_n\otimes A.$ Put
$A_n=\overbrace{A\otimes A\otimes\cdots \otimes A}^n.$ Denote
$j_1: C_1\to A\otimes A$ by $j_1(c\otimes a)=j(c)\otimes a$ for
all $a\in A$ and $c\in C.$ Define $j_{n+1}: C_{n+1}\to A_{n+2}$ by
$j_{n+1}(c\otimes a)=j_n(c)\otimes a$ for all $c\in C_n$ and $a\in
A,$ $n=1,2,....$

Note that, since $A$ is a unital simple AF-algebra,
$K_1(C_1)=K_1(C\otimes {\cal  U})\otimes K_0(A).$ Since
$K_1(C\otimes {\cal  U})$ and $K_0(A)$ are  torsion free, the map
$\phi_0: C\otimes {\cal  U}\to C\otimes {\cal  U}\otimes A$
defined by $a\to a\otimes 1$ induces an injective \hm\,
$(\phi_0)_{*1}.$

Define $\phi_1: C_1\to C_2$ by
\beq\label{kemb1}
\phi_1(a)={\rm diag}(a, 1\otimes j_1(a))\tforal a\in C_1\andeqn
\eneq
\beq\label{kemb2}
\phi_{n}(a)={\rm diag}(a, \overbrace{1\otimes j_n(a), 1\otimes
j_n(a),...,1\otimes j_n(a)}^{n-1})\tforal a\in C_n
\eneq
$n=1,2,....$ Note that, for each $n>1,$ $C_n\cong C_1\otimes A_n.$
$(\phi_n)_{*1}$ is injective, $n=1,2,....$
Put $B=\lim_{n\to\infty}(C_n, \phi_n).$ Then $\phi_{1, \infty}:
C_1\to B$ is a unital monomorphism. In particular,
$(\phi_{1, \infty})_{*1}|_{K_1(C_1)}$ is injective.
Therefore $(\phi_{1, \infty}\circ \phi_0)_{*1}$ is injective on
$K_1(C\otimes {\cal  U}).$  Since $A$ is a unital simple \CA, it is
also easy and standard to check that $B$ is a unital simple \CA.

To verify  that $B$ has tracial rank zero, fix any $\ep>0$ and a
finite subset ${\cal  F}\subset B$ and a positive element $a\in
B_+\setminus \{0\}.$ Without loss of generality, we may assume
that ${\cal  F}\subset \phi_{n,\infty}(C_n)$ for some large $n.$
There exists a finite subset ${\cal  G}\subset C_n$ such that
$\phi_{n, \infty}({\cal  G})\supset {\cal  F}.$ Fix an integer $k\ge
1.$ Then one may write
\beq\label{kemb5}
\phi_{n, n+k}(x)=\imath(x)\oplus \Phi(x)\rforal x\in C_n,
\eneq
where $\imath: C_n\to (1-p)C_{n+k}(1-p)$ is a unital monomorphism
such that $1-p=e_{11},$ where $\{e_{ij}\}$ is a system of matrix
unit for $M_{n(n+1)\cdots n(+k)}$ and $\Phi: C_n\to F,$ where $F$
is a unital AF-subalgebra of $pC_np$ (with $1_F=p$).

Let $F_{\infty}=\phi_{n+1,\infty}(F).$ Then $F_{\infty}$ is an
unital AF-algebra with $1_{F_{\infty}}=\phi_{n+1, \infty}(p).$ We
also have
\beq\label{kemb6}
[y, \, \phi_{n, \infty}(p)]=0\rforal y\in {\cal  F}.
\eneq
This, in particular, implies that $B$ satisfies the Popa condition
(1) and (2) as in 3.6.2 of \cite{Lnbk} (see also  \cite{Po}).
By 3.6.6 of \cite{Lnbk}, $B$ has the property
(SP). Now let $e\in {\overline{aBa}}$ be a non-zero projection. We
will show that we can choose $p$ (by choosing larger $k$ above) so
that
\beq\label{kemb7}
[1-p]\le [e].
\eneq
To do this, by replacing $e$ by an equivalent projection, we may
assume that $e\in \phi_{m, \infty}(C_m)$ for some $m\ge 1.$ Note
that each $\phi_n$ is injective. So there is a projection $q_0\in
C_m$ such that $\phi_{m, \infty}(q_0)=e.$ Note that
\beq\label{kemb8}
\phi_{m,m+1}(q_0)={\rm diag}(q_0, \overbrace{1\otimes
j_m(q_0),1\otimes j_m(q_0),...,1\otimes j_m(q_0)}^m)
\eneq
In particular,
\beq\label{kemb9}
[\phi_{m,\infty}(q_0)]\ge  [\phi_{m+1,\infty}(0,1\otimes
j_m(q_0),0,...,0)]
\eneq
In other words, there is a projection $q_1\in 1\otimes A\subset
C_{m+1}$ such that $[\phi_{m+1, \infty}(q_1)]\le [e].$ It then
easy to check that, for some sufficiently large $k,$ for the above
chosen $p,$
\beq\label{kemb10}
[1-p]\le [\phi_{m+1, \infty}(q_1)]\le [e].
\eneq
It follows that $B$ has tracial rank zero. Since each $C_n$
satisfies the UCT, so does $B$ (see \cite{RS}).

It is also standard to show that $B$ has a unique tracial state if $A$ does (see, for example, 3.7.10 of \cite{Lnbk}).

\end{proof}

\begin{cor}\label{Kemb2}
Let $C$ be a unital separable \CA\, in ${\cal  N},$ let $A\cong A\otimes {\cal  U}$ be a unital simple AF-algebra with a unique tracial state $\tau$ and
with $K_0(A)=\rho_A(K_0(A))$ which is  a subring of $\R$ with $[1_A]=1$ and
let $j: C\to A$ be a unital embedding.

Then there exists a unital separable simple \CA\, $B\in {\cal  N}$ with tracial rank zero,
with a unique tracial state, $\rho_B(K_0(B))=\rho_A(K_0(A))$ and a unital monomorphism
$\phi: C\to B$ satisfying the following:

\beq\label{Kemb2-1}
(\phi)_{*1}=\gamma\circ \kappa
\eneq
where $\kappa: K_1(C)\to K_1(C)/{\rm Tor}(K_1(C))$ is the quotient
map and $\gamma: K_1(C)/{\rm Tor}(K_1(C))\to K_1(B)$ is injective.
Moreover, there is a unital monomorphism $h: B\to A$ and  a unital
\hm\, with $h_{*0}=\rho_B$ such that
\beq\label{Kemb2-2}
(h\circ \phi)_{*0}=j_{*0}.
\eneq

\end{cor}

\begin{proof}
By considering the composition $C\to C\otimes {\cal  U}\to
A\otimes {\cal  U}\cong A,$  we may assume that $K_1(C)$ is
torsion free.
Since $C$ can be embedded into a unital simple AF-algebra, it has
a faithful tracial state. We will keep the notation used in the
proof of \ref{Kemb1}. Note that $C_{n+1}\cong C_n\otimes A$ for
each $n$ and for each $t\in T(C\otimes A),$ it has the form
$t_1\otimes \tau,$ where $t_1$ is a tracial state of $C$ and
$\tau\in T(A)$ is the unique tracial state. We will use $\tau$ for
the unique tracial state on $A_n.$

The \hm\, $\phi_n: C_n\to C_{n+1}$ now has the from
\beq\label{Kemb2-3}
\phi_n(a)={\rm diag}(a, \overbrace{1\otimes j_n(a), 1\otimes
j_n(a),...,1\otimes j_n(a)}^{n-1})
\eneq
for all $a\in C_n.$ It follows that
\beq\label{Kemb2-4}
(t_1\otimes
\tau)(\phi_n(a))={1\over{n}}t_1(a)+{n-1\over{n}}\tau(j_n(a))\rforal
a\in C_n .
\eneq
Thus, we compute that
\beq\label{Kemb2-5}
\lim_{n\to\infty}(t_1\otimes \tau)(\phi_{k,
n}(a))=\tau(j_k(a))\tforal a\in C_k.
\eneq
Since $B$ has a unique tracial state, we compute that
\beq\label{Kemb2-6}
\tau_1(\phi(a))=\tau(j(a))\tforal c\in C,
\eneq
where $\tau_1$ is the unique tracial state of $B.$ Since $B\in
{\cal N},$ by the classification theorem of \cite{Lnduke}, there
is a monomorphism from $B$ to $B_1$ which is a unital simple
AH-algebra with a unique tracial state, with
$K_0(B_1)=\rho_{B_1}(K_0(B_1))=\rho_B(K_0(B)).$ Then
(\ref{Kemb2-2}) follows.

\end{proof}

\begin{NN}\label{prebot}
{\rm Let $C_0$ be a unital amenable \CA, let $\Lambda: \Z^k\to Aut(C_0)$ ($k\ge 1$) be a \hm, and let $C=C_0\rtimes_{\Lambda}\Z^k.$
 Suppose  that $A=A\otimes {\cal U
}$ is a unital simple AF-algebra with a unique tracial state
$\tau$ and $K_0(A)=\rho_A(K_0(A))$ which is a subring of $\R$ with
$[1_A]=1$ and that $j: C\to A$ is  a unital monomorphism. By
\ref{Kemb2}, there exists a unital separable simple \CA\, $B_0\in
{\cal  N}$ with tracial rank zero, with a unique tracial state and with
$\rho_B(K_0(B_0))=\rho_A(K_0(A)),$ and there exists a unital monomorphism $j_1:
C\to B_0$ satisfying the following:

\beq\label{prebot-1}
(j_1)_{*1}=\gamma\circ \kappa
\eneq
where $\kappa: K_1(C)\to K_1(C)/{\rm Tor}(K_1(C))$ is the quotient
map and $\gamma: K_1(C)/{\rm Tor}(K_1(C))\to K_1(B_0)$ is injective.

Note also  in the proof of \ref{Kemb1} $K_1(C_n)=K_1(C)\otimes K_0(A)$ and the \hm\, $\phi_n: K_1(C_n)\to K_1(C_{n+1})$ is the identity
map. So $K_1(B_0)=K_1(C)\otimes K_0(A).$

Moreover, there is a unital monomorphism $h: B_0\to A$
with $h_{*0}=\rho_{B_0}$ such that
\beq\label{prebot-2}
(h\circ j_1)_{*0}=j_{*0}.
\eneq


Note that $K_0(A)$ is also divisible.
It follows that
\beq\label{prebot-3}
[h\circ j_1]=[j]\,\,\,\text{in}\,\,\, KK(C, A)\andeqn \tau(h\circ
j_1)=\tau\circ j.
\eneq
There exists (see \ref{drot} and \cite{KK2}) a \hm\, $R_{j, h\circ
\phi}: K_1(C)\to Aff(T(A))=\R.$ Then $R_{j, h\circ j_1}(K_1(C))$
is a countable subgroup of $\R.$ Fix a countable subring
$\D\subset \R$ which contains $\rho_A(K_0(A))$ and $R_{j, h\circ
j_1}(K_1(C))$ which also has the property that $\D\Q=\D.$  Let
$B=B(j, A)$ be a unital simple AF-algebra with unique tracial
state  and $K_0(B)=\rho_B(K_0(B))=\D$ and with $[1_B]=1.$  There
is a unital embedding $j_2: A\to B$ such that $(j_2)_{*0}={\rm
id}_{\rho_A(K_0(A))}.$ We may assume that $A\subset B.$

An important fact is that, if we regard $j: C\to B$ and $h\circ j_1: C\to B,$ then
\beq\label{prebot-4}
R_{j, h\circ j_1}(K_1(C))\subset \D=\rho_B(K_0(B)).
\eneq
}
\end{NN}

\begin{thm}\label{Bott1}
Let $C_0$ be a unital AH-algebra, $C=C_0\rtimes_{\af} \Z$ for some
$\af\in Aut(C_0),$  let $A\cong A\otimes {\cal  U}$ be a unital
simple AF-algebra with a unique tracial state and with
$K_0(A)=\rho_A(K_0(A))$ which is a subring of $\R$ and let $j:
C\to A$ be a unital monomorphism. Let $B=B(j, A)$ and $\D=K_0(B)$
be as in {\rm \ref{prebot}}.

Let $\ep>0,$ let ${\cal  F}\subset C$ be a finite subset and let ${\cal
P}\subset K_1(C).$ There is $\dt>0$  and a finite subset
$\{y_1,y_2,...,y_m\}\subset K_1(C)$ satisfying the following:

 If $\bt: K_1(C)\to K_0(A)$ is a \hm\, with
 \beq\label{Bott1-1}
 \rho_A((\bt(y_j))(\tau)<\dt,\,\,\,j=1,2,...,m
 \eneq
 where $\tau$ is the unique tracial state of $A,$
 then there exists a    unitary $u\in B$ such that
 \beq\label{Bott1-2}
 \|[\psi(a), \, u]\|<\ep\tforal a\in {\cal  F}\andeqn
 {\rm{bott}}_1(\psi,\, u)|_{{\cal  P}}=\bt'|_{{\cal  P}}
 \eneq
where $\psi: C\to B\otimes {\cal  U}\to B$ is defined by
$\psi(c)=j(c)\otimes 1_{\cal U}$ for $c\in C_0$ and $\psi(u_\af)=j(u_\af)\otimes
\imath(u_\sigma)$ and $\bt'=(j_2)_{*0}\circ \bt$ ($j_2$ is as in \ref{prebot}).
\end{thm}

 \begin{proof}
 Since $K_0(A)$ is torsion free, $\bt({\rm Tor}(K_1(C)))=\{0\}.$
 Therefore, there exists\\
  $\bt_1: K_1(C)/{\rm Tor}(K_1(C))\to K_0(A)$ such that
 \beq\label{Bot-1}
\bt=\bt_1\circ \kappa,
\eneq
where $\kappa: K_1(C)\to K_1(C)/{\rm Tor}(K_1(C))$ is the quotient
map.

Let $B_0$ be as in \ref{prebot}. By the classification of amenable
simple \CA s of tracial rank zero (\cite{Lnduke}), $B_0$ is a
unital simple AH-algebra of real rank zero. Write
$B_0=\lim_{n\to\infty}(C_n, \phi_n),$ where each $C_n$ has the
form $P_nM_{k(n)}(C(X_n))P_n$ and where $X_n$ is a finite CW
complex and $P_n\in M_{k(N)}(C(X_n))$ is a projection. Moreover,
let $\gamma$ and $j$ be as in \ref{prebot}.

Identify $K_0(A)$ with $\rho_A(K_0(A))\subset \R.$ We define $\gamma_1: K_1(B_0)=K_1(C)\otimes \rho_A(K_0(A))\to K_0(A)$ by
\beq\label{Bot-4-1}
\gamma_1(g\otimes r)=r\bt(g)\tforal r\in \rho_A(K_0(A))\andeqn g\in K_1(C).
\eneq
Then
\beq\label{Bot-4}
\gamma_1\circ \gamma=\bt_1.
\eneq

Furthermore, for any $x\in K_1(B_0),$ there are $g_1,g_2,...,g_l\in K_1(C)/{\rm Tor}(K_1(C))$ and real numbers
$r_1,r_2,...,r_l\in \rho_A(K_0(A))$ such that
\beq\label{Bot-4+}
\gamma_1(x)=\sum_{i=1}^l r_i\bt_1(g_i).
\eneq

Since $B_0$ satisfies the Universal Coefficient Theorem, there is a \hm\\ ${\tilde \bt}\in Hom_{\Lambda}(\underline{K}(B_0), \underline{K}(SB))$ such that
\beq\label{Bot-4+2}
{\tilde \bt}|_{K_1(B_0)}=\gamma_1.
\eneq

Without loss of generality, we may assume that $j_1({\cal  F})\subset \phi_{n, \infty}(C_n)$ and
$j_1({\cal  P})\subset [\phi_{n,\infty}](\underline{K}(C_n)).$

Let $\dt_0>0$ (in place of $\dt$) and $k(n)\ge n$ and $x_1, x_2,...,x_k$ be a set of generators for  $K_1(C_{k(n)})$
required by Lemma 7.5 of \cite{Lnasy} corresponding to $\ep/2$ and ${\cal  F}_1=h({\cal  F}).$
By (\ref{Bot-4+}), there are $y_1,y_2,...,y_m\in K_1(C)$ such that
\beq\label{Bot-4+3}
\gamma_1(x_j)=\sum_{i=1}^m r_{i,j}\bt(y_i),\,\,\,j=1,2,...,k.
\eneq
Therefore there exists $\dt>0$ such that
\beq\label{Bot-4+4}
\rho_B(\bt(y_j))(\tau)<\dt,\,\,\,j=1,2,...,m,
\eneq
implies
\beq\label{Bot-4+5}
\rho_B(\gamma_1(x_i))(\tau)<\dt_0,\,\,\,i=1,2,...,k.
\eneq
We now choose the above $\dt$ and $y_1,y_2,...,y_m$ and assume (\ref{Bott1-1}) holds.

Thus, by 7.5 of \cite{Lnasy}, there is a unitary $u_0\in
A$ such that
\beq\label{Bot-5}
\|[h, \, u_0]\|<\ep/2\tforal a\in {\cal  F}_1\andeqn \text{bott}_1(h, u_0)|_{\cal  P}={\bar \bt}|_{\cal  P}.
\eneq
Viewing  $A$ as a unital \SCA\, of $B,$
 by \ref{prebot},
\beq\label{Bot-6}
[h\circ j_1]=[j]\,\,\,\text{in}\,\,\, KK(C, B),\,\,\,\tau\circ
h\circ j_1=\tau\circ j\andeqn
\eneq
\beq\label{Bot-7}
R_{j, h\circ j_1}(K_1(C))\subset \rho_B(K_0(B)).
\eneq

Thus, by \ref{2Tuni}, there is a unitary $w\in U(B)$ such that
\beq\label{Bot-8}
{\rm ad}\, w\circ \psi_1\approx_{\ep/2} \psi\,\,\,\text{on}\,\,\,
{\cal  F},
\eneq
where $\psi_1(c)=h\circ j_1(c)\otimes 1$ and $\psi_1(u_\af)=h\circ
j_1(u_\af)\otimes \imath(u_\sigma).$
 Choose $u=w^*u_0w.$ Then
\beq\label{Bot-9}
\|[\psi(a), \, u]\|<\ep\tforal a\in {\cal  F}\andeqn
\text{bott}_1(\psi, u)|_{{
\cal P}}=\bt|_{{\cal  P}}.
\eneq

\end{proof}

\section{Homotopy lemmas}

The next follows immediately from Theorem 3.9 of \cite{Lnuct}.
Several versions of the following have been appeared.

\begin{thm}\label{QUNI}
Let $C$ be a unital separable amenable \CA. For any $\ep>0$ and any
finite subset ${\cal  F}\subset C,$ there exists $\dt>0,$ a finite subset ${\cal  G}\subset C,$
a finite subset ${\cal  Q}\subset \underline{K}(C)$ and an integer $n\ge 1$ satisfying the following:

For any unital simple \CA\, $A$ of tracial rank zero and any two unital $\dt$-${\cal  G}$-multiplicative
\morp s $L_1, L_2: C\to A$ and any
unital monomorphism $h: C\to A,$
if
\beq\label{quni1}
[L_1]|_{\cal  Q}=[L_2]|_{\cal  Q},
\eneq
then there exists a unitary $U\in M_{n+1}(A)$ such that
\beq\label{quni2}
{\rm ad}\, U\circ {\bar L}_1\approx_{\ep} {\bar L}_2\,\,\,\text{on}\,\,\, {\cal  F},
\eneq
where
$$
{\bar L}_i(a)={\rm diag}(L_1(a), \overbrace{h(a),h(a),...,h(a)}^n)
$$ for all $a\in C$ and $i=1,2.$
\end{thm}

\begin{proof}
We will apply Theorem 1.2 of \cite{GL1}. By Theorem 3.9 of
\cite{Lnuct} and the assumption that there exists a unital
embedding $h: C\to A,$  $C$ satisfies property (P) of Remark 1.1
of \cite{GL1}. Thus Theorem 1.2 of \cite{GL1} holds for $C$ (in
place $A$) and $A$  (in place of $B$). Since in this case $A$ has
tracial rank zero,  weakly unperforated $K_0(A),$  stable rank one
and real rank zero, the theorem follows  ( see also the proof of
Theorem 3.1 of \cite{GL1}).

\end{proof}

\begin{lem}\label{Key}
Let $C_0$ be a unital AH-algebra, $C=C_0\rtimes_\af\Z$ for some
$\af\in Aut(C_0)$  and let $A\cong A\otimes {\cal U}$ be a
unital simple AF-algebra with a unique tracial state $\tau$ and
with $K_0(A)=\rho_A(K_0(A)).$ Let  $j: C\to A$ be a unital
monomorphism, $\ep>0,$ ${\cal  F}\subset C_0$ be a finite subset
and ${\cal  P}\subset K_1(C_0)$ be a finite subset.
 There
is $\dt>0,$ a finite subset ${\cal  G}\subset C_0,$ and integer $K$ and a finite
subset ${\cal  Q}\subset \underline{K}(C)$ satisfying the
following:

Suppose that $v\in U(A)$ is a unitary such that
\beq\label{key1}
\|[j(a),\, v]\|<\dt\tforal a\in {\cal  G}\andeqn
\eneq
\beq\label{key2}
[L]|_{\cal  Q}=[j]|_{\cal  Q},
\eneq
where $L: C\to A$ is a \morp\, such that
\beq\label{key3}
\|L(\sum_{k=-K}^{K} f_ku_\af^k)-\sum_{k=-K}^{K}j(f_k)(j(u_\af)v)^k\|<\dt
\eneq
for all $f_k\in {\cal G}.$
Then there exists a unitary $w\in U(A)$ such that
\beq\label{key4}
\|[j(a),\, w]\|<\ep\tforal a\in {\cal F}\andeqn
{\rm{bott}}_1(j|_{C_0}, j(u_\af)^*w^*j(u_\af)vw)|_{\cal  P}=0.
\eneq

\end{lem}

\begin{proof}

Let $\ep>\ep_0>0$ and ${\cal  F}\subset {\cal  F}_1\subset C_0$ be a
finite subset. We assume that
\beq\label{key-1}
\text{bott}_1(h,\, v')|_{\cal  P}
\eneq
is well defined for any unital \hm\, $h$ from $C_0$ and any
unitary $v'$ provided that $\|[h(a), \, v']\|<\ep_0$ for all $a\in
{\cal  F}_1.$


Let ${\cal P}=\{x_1,x_2,...,x_k\}.$

We write $C_0=\lim_{n\to\infty}(C_n, \phi_n),$ where each $C_n$
has the form $P_nM_{r(n)}(C(X_n))P_n,$ where $X_n$ is a finite CW
complex and $P_n\in M_{r(n)}(C(X_n))$ is a projection.  We may
assume that ${\cal  P}\subset [\phi_{n, \infty}](K_1(C_n)).$ Let
$\eta>0$ (in place of $\dt$), $k(n)\ge n$ and $y_1,y_2,...,y_m\in
K_1(C_{k(n)})$ be a set of generators as required by Lemma 7.3 of
\cite{Lnasy} for $\{x_1,x_2,...,x_k\},$ $\ep_0/2$ and ${\cal
F}_1.$ Denote $z_j=(\phi_{k(n), \infty})_{*1}(y_j),$
$j=1,2,...,m.$ We may assume that ${\cal P}\subset
\{z_1,z_2,...,z_m\}.$

Let ${\cal  P}_1'\subset \underline{K}(C_0)$ be a finite subset
containing $z_1,z_2,...,z_m$
and ${\cal
P}_1={\cal  P}_1'\cup [\af]({\cal  P}_1').$

Let $\ep_1>0$ and ${\cal  F}_2\subset C_0$ be a finite subset such
that
\beq\label{key-2}
\text{bott}_1(h',\, v')|_{{\cal  P}_1}
\eneq
is well defined for any unital \hm\, $h'$ from $C_0$ and any
unitary $v'$ provided that $\|[h'(a), \, v']\|<\ep_1$ for all
$a\in {\cal  F}_2.$ We may also
\beq\label{key-3}
\text{bott}_1(h'', v')|_{{\cal  P}_1}=\text{bott}_1(h',v')|_{{\cal
P}_1}=\text{Bott}(h', v'')|_{{\cal  P}_1},
\eneq
provided that $h''\approx_{\ep_1} h'$ on ${\cal  F}_2$ and $\|v' -v''\|<\ep_1.$ We may
further assume that $\ep_1<\ep_0/2$ and ${\cal  F}_1\subset {\cal
F}_2.$

 To simplify notation, we may also assume that
\beq\label{key-4}
\rho_D(\text{bott}_1(h', v'))(z_j)(t)<\eta,\,\,\,j=1,2,...,m
\eneq
for any unital \hm\, $h': C\to D,$ any unitary $v'\in D$ and
unital \CA\, $D$ and tracial state $t\in T(D),$ provided that
\beq\label{key-5}
\|[h'(a),\, v']\|<\ep_1\rforal a\in {\cal  F}_2.
\eneq
Put ${\cal  G}_1={\cal  F}_2\cup\{u_\af\}.$

 Let $\dt_1>0$ ( in place of $\dt$)  a finite subset
${\cal  G}_2\subset C$ (in place of ${\cal  G}$) be a finite subset, ${\cal  Q}\subset
\underline{K}(C)$
 be a finite subset and $N\ge 1$ be an integer required by \ref{QUNI}
 for $\ep_1/2$ and ${\cal  G}_1.$

There is $\dt,$  a finite subset ${\cal  G}\subset  C$ and an integer $K\ge 1$ such that $L$ is $\dt_1$-${\cal  G}_1$-multiplicative
if (\ref{key1}) and (\ref{key3}) hold.

Suppose (\ref{key1}), (\ref{key2}) and (\ref{key3}) hold.

It follows from \ref{QUNI}  that there exists a unitary $U\in M_{N+1}(A)$
such that
\beq\label{key-6}
{\rm ad}\, U\circ {\bar L}\approx_{\ep_1/2} {\bar
j}\,\,\,\text{on}\,\,\, {\cal  G}_1,
\eneq
where
\beq\label{key-7}
{\bar L}(a)&=&{\rm diag}(L(a),
\overbrace{j(a),j(a),...,j(a)}^N)\andeqn\\
 {\bar j}(a)&=& {\rm diag}(\overbrace{j(a),j(a),...,j(a)}^{N+1})
 \eneq
 for all $a\in C.$
It follows that
\beq\label{key-8}
\|[{\bar j}(a),\, U]\|<\ep_1/2\tforal a\in {\cal  F}_2.
\eneq

Put $V={\rm diag}(\overbrace{v,v,...,v}^{N+1}).$ By (\ref{key-6}),
\beq\label{key-9}
\|1-{\bar j}(u_\af)^*U^*{\bar j}(u_\af)VU\|<\ep_1/2.
\eneq
It follows that
\beq\label{key-10}
\text{bott}({\bar j}|_{C_0},\,{\bar j}(u_\af)^*U^*{\bar
j}(u_\af)VU)|_{{\cal  P}_1}=0.
\eneq
Thus (by (\ref{key-3})),
\beq\label{key-11}
\hspace{-0.2in}0\hspace{0.2in}&\hspace{-0.3in}=&\text{bott}_1({\bar
j}|_{C_0},\, V)|_{{\cal  P}_1}+ \text{bott}_1({\bar j}|_{C_0},\,
V^*{\bar j}(u_\af)^*U^*{\bar
j}(u_\af)VU)\\
&\hspace{-0.4in}=&\text{bott}_1({\bar j}|_{C_0},\, V)|_{{\cal
P}_1} +\text{bott}_1({\bar j}|_{C_0},\,{\bar j}(u_\af)^*U^*{\bar
j}(u_\af)VUV^*)|_{{\cal  P}_1}\\
&\hspace{-0.4in}=& \text{bott}_1({\bar j}|_{C_0}, V)|_{{\cal
P}_1}+ \text{bott}_1({\bar j}|_{C_0},\,{\bar j}(u_\af)^*U^*{\bar
j}(u_\af))|_{{\cal  P}_1}+\text{bott}_1({\bar
j}|_{C_0},VUV^*)|_{{\cal  P}_1}\\\label{key-12}
&\hspace{-0.4in}=&\text{bott}_1({\bar j}|_{C_0}, V)|_{{\cal P}_1}+
\text{bott}_1({\bar j\circ \af^{-1}}|_{C_0},\,U^*)|_{{\cal
P}_1}+\text{bott}_1({\bar j}|_{C_0},U)|_{{\cal  P}_1}
\eneq

Let $G_1$ be a subgroup of $K_1(C)$ generated by ${\cal  P}_1.$
Since $K_0(A)$ is divisible, there are \hm s $\lambda_1, \gamma:
K_1(C_0)\to K_0(A)$ which extend $\text{bott}_1({\bar j}|_{C_0},\,
V)$ and $\text{bott}_1({\bar j}|_{C_0},\, U),$ respectively.

It follows from (\ref{key-12}) that
\beq\label{key-13}
\lambda_1|_{G_1}=-\gamma|_{G_1}+\gamma\circ \af^{-1}_{*1}|_{G_1}.
\eneq
It should be note that
\beq\label{key-14}
(\lambda_1)|_{G_1}=(N+1)\text{bott}_1(j|_{C_0}, v)|_{G_1}.
\eneq

 Define $\gamma_1: K_1(C_0)\to K_0(A)$ by
$\gamma_1(x)=-{1\over{N+1}}\gamma(x)$ for all $x\in K_1(C_0).$

By (\ref{key-4}) and applying Lemma 7.3 of \cite{Lnasy}, there is
a unitary $w\in A$ such that
\beq\label{key-14+}
\|[j(a), \, w]\|<\ep_0/2\andeqn \text{bott}_1(j,\,
w)(x_j)=\gamma_1(x_j),\,\,\,j=1,2,...,k.
\eneq

We then compute (using (\ref{key-13}) among other things) that
\beq\label{key-15}
&&\hspace{-0.6in}\text{bott}_1(j,\, j(u_\af^*)w^*j(u_\af)vw)(x_j)\\
&=&\text{bott}_1(j,v)(x_j)+\text{bott}_1(j,\, v^*j(u_\af)^*w^*
j(u_\af)vw)(x_j)\\
&=&\text{bott}_1(j,\, v)(x_j)
+\text{bott}_1(j,\,j(u_\af)^*w^*j(u_\af)vwv^*)(x_j)\\
&=& \text{bott}_1(j, v)(x_j)+
\text{bott}_1(j,\,j(u_\af)^*w^*j(u_\af))(x_j)+\text{bott}_1(j
,wvw^*)(x_j)\\\label{key-16}
&=&\text{bott}_1(j, v)(x_j)+
\text{bott}_1(j\circ \af^{-1},\,w^*)(x_j)+\text{bott}_1(
j,w)(x_j)\\
&=&{1\over{N+1}}(\lambda_1)(x_j)+(-\gamma_1\circ
(\af^{-1}_{*1}(x_j)+\gamma_1(x_j))\\
&=&({1\over{N+1}})(\lambda_1-\gamma+\gamma\circ
\af^{-1}_{*1})(x_j)=0,
\eneq
$j=1,2,...,k.$

It follows that
\beq\label{key-17}
\text{bott}_1(j_{C_0}, \, j(u_\af^*)w^*j(u_\af)vw)|_{\cal  P}=0.
\eneq

\end{proof}

\begin{thm}\label{S1uni}
Let $C$ be a unital AH-algebra
and let $\af\in Aut(C)$ be an
automorphism and let $A\cong A\otimes {\cal U}$ be a unital simple AF-algebra with a unique
tracial state $\tau$ and $\rho_A(K_0(A))=K_0(A).$
Let $h: (C\rtimes_\af\Z)\otimes C(\T)\to A$ be a unital monomorphism.

For any $\ep>0$ and a finite subset ${\cal  F}\subset C\rtimes_\af
\Z,$ there exists $\dt>0,$  $\eta>0,$ a finite subset ${\cal  G}\subset
C\rtimes_\af \Z,$  a finite subset ${\bar {\cal G}}\subset C\otimes C(\T)$ and a finite subset ${\cal  P}\subset
K_1(C\rtimes_\af \Z)$ satisfying the following:

Suppose that there is a unitary $v\in U(A)$ and  \morp\, $L: C\otimes C(\T)\to A$  such
that
\beq\label{S1u1}
&&\|[h(a), \, v]\|<\dt\tforal a\in {\cal  G},
{\rm{bott}}_1(h, v)|_{\cal  P}=0\andeqn\\\label{S1u1+}
&&\tau\circ h|_{C\otimes C(\T)}\approx_{\eta}\tau \circ L\,\,\,\text{on}\,\,\,{\bar {\cal G}},
\eneq
and
\beq\label{S1u2}
L\approx_{\dt} h\,\,\,\text{on}\,\,\, {\cal G}\andeqn L(1\otimes z)\approx_{\dt} v
\eneq
and where $z\in C(\T)$ is the identity function on the unit circle.
Then there exists a unitary $W\in A\otimes{\cal U}$ such that
\beq\label{S1u3}
\|[h^{(1)}(c), \,W]\|<\ep\rforal c\in {\cal F}\andeqn W^*(v\otimes 1_{\cal U})W\approx_{\ep} h(1\otimes z)\otimes 1_{\cal U},
\eneq
where $h^{(1)}(c)=h(c)\otimes 1_{\cal U}$ for all $c\in C$ and $h^{(1)}(u_\af)=h(u_\af)\otimes \imath(u_\sigma).$
\end{thm}

\begin{proof}

Let $\ep>0$ and ${\cal  F}\subset C\rtimes_\af \Z.$ Without loss
of generality, we may assume that ${\cal  F}={\cal
F}_0\cup\{u_\af\}$ for some finite subset ${\cal  F}_0\subset C.$
Define $C_1=C\otimes C(\T).$ Define $\af'\in Aut(C_1)$ by
$\af'(c\otimes 1_{C(\T)})=\af(c)\otimes 1_{C(\T)}$ and
$\af'(1_C\otimes f)=1_C\otimes f$ for all $f\in C(\T).$ In other
words, $C_1\rtimes_{\af'} \Z \cong (C\rtimes_\af\Z)\otimes C(\T).$
Define ${\cal  F}'={\cal  F}_0\cup\{1\otimes z\}.$

To apply \ref{2burg}, let $\dt_1>0$ (in place of $\dt$), ${\cal
G}_1\subset C_1$ (in place of ${\cal  G}$) be a finite subset and
${\cal  P}_1'\subset \underline{K}(C_1)$ (in place of ${\cal  P}$)
be a finite subset required by \ref{2burg} for $\ep/2$ and ${\cal
F}'.$ We may assume that $\dt_1<\ep/2$ and ${\cal  G}_1\supset {\cal
F}'.$ Since $K_1(A)=\{0\}$ and $K_0(A)$ is torsion free and
divisible, we may assume that ${\cal P}_1'={\boldsymbol{\bt}}({\cal
P}_1)$ for some finite subset ${\cal P}_1\subset K_1(C_1).$

Let $\dt_2>0$ (in place of $\dt$), ${\cal  G}_2\subset C_1$ be a
finite subset, $K\ge 1$ be an integer and ${\cal  Q}\subset
\underline{K}(C_1\rtimes_{\af'}\Z)$ be a finite subset required by
\ref{Key} for $\dt_1/2,$ ${\cal  P}_1$ and ${\cal  G}_1.$  Since
$K_0(A)$ is torsion free and divisible and $K_1(A)=\{0\},$ we may
assume that ${\cal  Q}\subset
K_0(C_1\rtimes_{\af'}\Z)=K_0((C\rtimes_\af\Z)\otimes C(\T)).$ We
may assume that ${\cal  Q}={\cal  P}_2\cup{\boldsymbol{\bt}}
({\cal  P}_3),$ where ${\cal  P}_2\subset K_0(C\rtimes_{\af}\Z)$
and ${\cal  P}_3\subset K_1(C\rtimes_\af\Z)$ are finite subsets.
We may assume that ${\cal G}_2\supset {\cal G}_1.$

By \ref{APP}, choosing $0<\dt_2<\dt_2'$ and a finite subset ${\cal
G}_2'\supset {\cal  G}_2$ so that there exits a unital
\morp\, $L_1: C_1\rtimes_{\af'}\Z \to A\otimes {\cal  U}$ such that
\beq\label{lh-2-1}
\|L_1(\sum_{k=-K}^Kg_ku_{\af}^k)-\sum_{k=-K}^Kh(g_k)(h(u_\af)v')^k\|<\dt_2/2
\eneq
for all $g_k\in {\cal G}_2$ and
for any $v'$ such that $\|[h(a), \, v']\|<\dt_2'$ for all $a\in
{\cal  G}_2'.$ We may assume that ${\cal  G}_2'={\cal
G}_2''\otimes \{1_{C(\T)},z\}$ for some finite subset ${\cal
G}_2''\subset C.$ Without loss of generality, we may also assume
that
\beq\label{lh-2-1+}
[L_1]|_{\cal  Q}=[L_1']|_{\cal  Q}
\eneq
provided that $L_1, L_1'$  both satisfy(\ref{lh-2-1})
$L_1\approx_{\dt_2} L_1'$ on ${\cal  G}_2\cup\{u_\af\}.$

 Define  $\tau_0=\tau\circ h$ on
$C_1.$

Let $0<\dt_3<\dt_2,$ $\eta>0$ (in place of $\sigma$), ${\bar {\cal
G}} \subset C_1$ (in place of ${\cal G}$) be a finite subset and
${\cal P}_4\subset \underline{K}(C_1)$ be a finite subset required
by \ref{2QUni} for $\min\{\dt_2/4, \ep/4)\}$ and ${\cal  G}_2'$ (with $\tau_0$ above).
Since $K_0(A)$ is torsion free and divisible and $K_1(A)=0,$ we
may assume that ${\cal  P}_4\subset  K_0(C_1).$ Without loss of
generality, (by choosing smaller $\dt_3,$ for example), we may
assume that ${\bar {\cal  G}}={\cal  G}_3'\otimes {\cal G}_3'',$
where $\{1_C\}\subset {\cal  G}_3'\subset C$ and
$\{1_{C(\T)},\,z\}\subset {\cal G}_3''\subset C(\T)$ are  finite
subsets, and we may assume that ${\cal  P}_4={\cal
P}_5\cup{\boldsymbol{\bt}}( {\cal  P}_6),$ where ${\cal
P}_5\subset K_0(C)$ and ${\cal  P}_6\subset K_1(C)$ are  finite
subsets.

Without loss of generality, we may assume that
\beq\label{lh-2++}
[L]|_{{\cal  P}_5}=[L']|_{{\cal  P}_5}
\eneq
provided that $L, L'$ are $\dt_3$-${\cal G}_3'$-multiplicative
\morp s from $C_1$ and $L\approx_{\dt_3} L'$ on ${\cal  G}_3'.$

Denote by $j_0: C\to C\rtimes_{\af}\Z$  the natural embedding.
Put ${\cal  P}=(j_0)_{*1}({\cal  P}_6)\cup {\cal  P}_3.$ Let $\eta_1>0$ such that $\eta_1<\min\{\dt_3/4,\dt_2/4,
\ep/2,\dt_1/4,\eta\}.$

It follows from \ref{APP} that there is $\dt_4>0,$  a finite subset ${\cal
G}_4\subset C$ and a unital  $\eta_1$-${\bar {\cal
G}}$-multiplicative \morp\, $L': C_1\to A$ such that
\beq\label{lh-3}
L'(c\otimes f)\approx_{\eta_1} h(c)f(v)\rforal c\in {\cal  G}_2'\andeqn  f\in {\cal G}_2'',
\eneq
provided that $\|[h(a), \, v]\|<\dt_4$ for all $a\in {\cal
G}_4.$ We may assume that ${\cal  G}_4\supset {\cal G}_2''\cup
{\cal  G}_2'''.$

Put ${\cal G}={\cal G}_4\cup\{u_\af\}.$

Put $\dt=\min\{\dt_4/2,\eta_1/2\}.$ Suppose that (\ref{S1u1}), (\ref{S1u1+}) and (\ref{S1u2})
hold for the above $\dt,$ $\eta,$ ${\cal  G},$ ${\bar {\cal G}}$ and ${\cal  P}$ (and for some $v$ and $L$).
Therefore, we may assume that $L$  is an $\eta_1$-${\cal  G}_2$-multiplicative \morp\,  satisfying (\ref{lh-3})
(replacing $L'$ by $L$).
We still have  that
\beq\label{lh-4}
\tau\circ L\approx_{\eta} \tau_0\,\,\,\text{on}\,\,\,{\bar {\cal  G}}.
\eneq


Note that
\beq\label{lh-7}
K_0(C_1)=K_0(C)\oplus K_1(C).
\eneq
We have (by (\ref{lh-2++}))
\beq\label{lh-8}
[L]|_{{\cal  P}_5}=h_{*0}|_{{\cal  P}_5}.
\eneq
By (\ref{lh-2++}) and (\ref{S1u1}) (and the choice of ${\cal P}$), we may also assume
that
\beq\label{lh-9}
[L]|_{{\boldsymbol{\bt}}({\cal  P}_6)}=0.
\eneq
Since ${\rm ker}\rho_A=\{0\},$
\beq\label{lh-10}
h_{*0}|_{{\boldsymbol{\bt}}(K_1(C))}=0.
\eneq
Therefore
\beq\label{lh-11}
[h]|_{{\cal  P}_4}=[L]|_{{\cal  P}_4}.
\eneq
Moreover,  by  (\ref{lh-4}),
\beq\label{lh-12}
\tau\circ h=\tau_0\approx_{\eta} \tau\circ
L\,\,\,\text{on}\,\,\, {\bar {\cal  G}}.
\eneq
By applying \ref{2QUni}, there exists a unitary $w_0\in U(A )$ such that
\beq\label{lh-13}
{\rm ad}\, w_0\circ L \approx_{\min\{\dt_2/4, \ep/2)\}}
h\,\,\,\text{on}\,\,\, {\cal  G}_2'.
\eneq
In particular,
\beq\label{lh-14}
w_0^*vw_0\approx_{\min\{\dt_2/4, \ep/2\}} h(1\otimes z)\andeqn \|[h(c),\,w_0]\|<\min\{\dt_2/4,\ep/2\}\rforal c\in {\cal G}_2''.
\eneq

Define $V=h(u_\af)^*w_0^*h(u_\af)w_0.$
Then
\beq\label{lh-15}
\|[h(a), \, V]\|<\min\{\ep/2, \dt_2/4\}\rforal a\in {\cal  G}_2'.
\eneq
By (\ref{lh-15}), there exists a unital  \morp\, $L_1:
C_1\rtimes_{\af'}\Z\to A\otimes {\cal U} $ which satisfies
(\ref{lh-2-1}) with $v'=V.$ Therefore
\beq\label{lh-16}
L_1\approx_{\dt_2} h\,\,\,\text{on}\,\,\,{\cal  G}_2''\andeqn
L_1(u_\af)\approx_{\dt_2} h(u_\af)V=w_0^*h(u_\af)w_0.
\eneq
Note that
\beq\label{lh-17}
[L_1]|_{{\cal  P}_2}=[h]|_{{\cal  P}_2}.
\eneq
By the assumption (\ref{S1u1}) and by   (\ref{lh-14}),
we have
\beq\label{lh-18}
\text{bott}_1({\rm ad}\, w_0\circ h,\, h(1\otimes z))|_{{\cal
P}_3}&= &\text{bott}_1({\rm ad}\, w_0\circ h,\, w_0^*vw_0)|_{{\cal
P}_3}\\\label{lh-18+}
 &=&\text{bott}_1(h,\, v)|_{{\cal  P}_3}=0.
\eneq
It follows from (\ref{lh-14}), (\ref{lh-16}), (\ref{lh-2-1+}) and
(\ref{lh-18+}) that
\beq\label{lh-19-1}
[L_1]|_{{\boldsymbol{\bt}}({\cal  P}_3)}=0.
\eneq
Since ${\rm ker}\rho_{A}(A )=\{0\},$ we have
\beq\label{lh-19}
[L_1]|_{{\boldsymbol{\bt}}({\cal  P}_3)}=0=[h]|_{{\boldsymbol{\bt}}({\cal  P}_3)}.
\eneq
Thus
\beq\label{lh-20}
[L_1]|_{\cal  Q}=[h]|_{\cal  Q}.
\eneq
Therefore, by \ref{Key}, there is a unitary $w_1\in U(A\otimes{\cal U} )$ such that
\beq\label{lh-21}
\|[h(a), \, w_1]|\|<\dt_1/2\rforal a\in {\cal  G}_1\andeqn \text{Bott}(h,\, h(u_\af)^*w_1^*h(u_\af)Vw_1)|_{{\cal  P}_1}=0.
\eneq
Note that $h(u_\af)V=w_0^*h(u_\af)w_0.$
It follows from \ref{2burg} that there exists $w_2\in U(A \otimes {\cal U})$ such that
\beq\label{lh-22}
\hspace{-1.2in} w_2^* (h(a)\otimes 1_{\cal U})w_2&\approx_{\ep/2}&h(a)\otimes 1_{\cal U}\rforal a\in
{\cal  F}'\andeqn\\\label{lh-22+}
 w_2^*(w_1^*(w_0^*h(u_\af)w_0)w_1)\otimes
\imath(u_\sigma))w_2&\approx_{\ep/2}& h(u_\af) \otimes
\imath(u_\sigma).
\eneq
Define $W=(w_0w_1\otimes 1_{\cal U})w_2.$
Then, by (\ref{lh-14}), (\ref{lh-22}) and (\ref{lh-22+}),
\beq\label{lh-23}
\|[h^{(1)}(c),\, W]\|<\ep\rforal c\in {\cal F}\andeqn W^*(v\otimes 1_{\cal U})W\approx_{\ep}h(1\otimes z)\otimes 1_{\cal U}.
\eneq

\end{proof}

The following holds in great generality and well-known. We only use
the special case below. The proof in this case is easier than the
general case.

\begin{lem}\label{meas}
Let $A$ be a unital simple AF-algebra with a unique tracial state
$\tau$ and $\rho_A(K_0(A))$ is a countable dense divisible subgroup
of $\R.$ There is a unitary $u\in A$ such that
\beq\label{meas1}
\tau\circ f(u)=\int_{\T} f(t) dt\rforal f\in C(\T),
\eneq
where the integral is  on the standard Lesbegue measure.

\end{lem}

\begin{lem}\label{2meas}
Let $A$ be a unital simple AF-algebra with a tracial state $t$ and
let $v\in U(A).$ Suppose that $u\in {\cal  U}$ is a unitary such
that
\beq\label{2meas1}
\tau\circ f(u)=\int_{\T} f(t) dt\rforal f\in C(\T).
\eneq
Suppose $C$ is a unital separable \CA\, and $h: C\to A$ is a unital
monomorphism. Define  a linear map $\phi: C\otimes C(\T)\to A\otimes
{\cal  U}$ by $\phi(c\otimes f)=h(c)f(v\otimes u)$ for all $c\in C$
and $f\in C(\T).$

 Then
\beq\label{2meas2}
(t\otimes \tau)(\phi(c\otimes f))=t(h(c))\cdot \int_\T f(t)dt\tforal
c\in C\andeqn f\in C(\T),
\eneq
where $\tau$ is the unique tracial state on ${\cal  U}.$

\end{lem}

\begin{proof}
Let $\ep>0$ and ${\cal  F}\subset C(\T)$ be a finite subset. We
assume that $z\in {\cal  F},$ where $z$ denotes the identity
function on the unit circle.

By \cite{LnFU}, for any $\dt>0,$ there are (for some large $n$),
mutually orthogonal projections $p_1,p_2,...,p_n$ in $A$ and
mutually orthogonal projections $e_1,e_2,...,e_n$ in ${\cal  U}$
with $\sum_{i=1}^np_i=1_A$ and $ \sum_{i=1}^n e_i=1_{\cal U}$ such
that
\beq\label{2meas-1}
f(v) &\approx_{\dt/4}& \sum_{k=1}^n f(\omega_k) p_k\andeqn\\
f(u)&\approx_{\dt/4}& \sum_{k=1}^n f(\omega_k) e_k
\eneq
for all $f\in {\cal  F},$ where $\omega_k=e^{2k\pi i\over{n}}$
 and $\tau(e_k)=1/n,$ $k=1,2,...,n.$
We have that
\beq\label{2meas-2}
(\sum_{k=1}^n\omega_k p_k)\otimes (\sum_{k=1}^n \omega_k e_k)=
\sum_{k=1}^n \omega_k(\sum_{j+i=k(\text{mod}\, n)} p_j\otimes e_i).
\eneq
Put $q_k=\sum_{j+i=k(\text{mod}\, n)} p_j\otimes e_i.$ Then
\beq\label{2meas-3}
(t\otimes \tau)(q_k)&=&\sum_{j+i=k(\text{mod}\, n)} t(p_j)\tau(
e_i)\\
&=&\sum_{j+i=k(\text{mod}\, n)} t(p_j)(1/n)=1/n.
\eneq
It follows immediately that, with sufficiently small $\dt$ (or
sufficiently large $n$),
\beq\label{2meas-4}
\tau (\phi(1\otimes f))\approx_{\ep} \tau(f(u))\rforal f\in {\cal
F}.
\eneq
Furthermore, assuming $\|c\|\le 1,$
\beq\label{2meas-5}
(t\otimes \tau)(\phi(c\otimes f))&\approx_{\ep}& (t\otimes
\tau)(\sum_{k=1}^n f(\omega_k)(h(c) q_k)\\
&=& \sum_{k=1}^n f(\omega_k) (t\otimes \tau)((h(c)q_k)\\
&=&\sum_{k=1}^n f(\omega_k)(1/n) t(h(c)\sum_{j=1}^n p_j)\\
&=&\sum_{k=1}^n f(\omega_k)(1/n)t(h(c))\\
&\approx_{\ep}& t(h(c))\tau(f(u)).
\eneq
for all $f\in {\cal  F}.$ Let $\ep\to 0,$ we conclude that
\beq\label{2meas-6}
(t\otimes \tau)(\phi(c\otimes f))=t(c)\tau(f(u))\tforal c\in
C\andeqn f\in C(\T).
\eneq

\end{proof}

\begin{lem}\label{Lhom}
Let $C$ be a unital AH-algebra
and let $\af\in Aut(A)$ be an automorphism and let $A\cong A\otimes
{\cal U}$ be a unital simple AF-algebra with a unique tracial state
$\tau$ and $\rho_A(K_0(A))=K_0(A).$ Let $\phi: C\rtimes_\af\Z\to A$
be a unital monomorphism.

For any $\ep>0$ and a finite subset ${\cal  F}\subset C\rtimes_\af
\Z,$ there exists $\dt>0,$ a finite subset ${\cal  G}\subset
C\rtimes_\af \Z$ and a finite subset ${\cal  P}\subset
K_1(C\rtimes_\af \Z)$ satisfying the following:

Suppose that there is a unitary $v\in U(A)$ such
that
\beq\label{Lh1}
\|[\phi(a), \, v]\|<\dt\tforal a\in {\cal  G}\andeqn
{\rm{bott}}_1(\phi, v)|_{\cal  P}=0.
\eneq

Then there exists a continuous path of unitaries $\{v_t: t\in [0,
1]\}\subset U(A\otimes  {\cal U})$ such that
\beq\label{Lh2}
v_0=v\otimes 1_{ {\cal U}},\,\,\, v_1=1_{A\otimes {\cal U}},\,\,\,\|[\psi(a), v_t]\|<\ep
\eneq
for all $a\in {\cal  F}$ and $t\in [0,1],$ and
\beq\label{Lh3}
{\rm{Length}}(\{v_t\})\le 2\pi+\min\{\ep,1\},
\eneq
 where $\psi:
C\rtimes_\af\Z\to A\otimes {\cal U}$ defined by $\psi(a)=a\otimes 1$ for
$a\in C$ and $\psi(u_\af)=\phi(u_\af)\otimes \imath(u_\sigma).$

\end{lem}

\begin{proof}
To simplify notation, without loss of generality, it suffices
to prove the theorem with $v_t\in A\otimes {\cal  U}\otimes  {\cal U}$ and
$\psi(u_\af)=\phi(u_\af)\otimes 1_{\cal  U}\otimes\imath(u_\sigma).$

 Let $u_0\in {\cal  U}$ be a unitary with
\beq\label{lh-1}
\tau_1(f(u_0))=\int_{\T} f(t) dt\tforal f\in C(\T),
\eneq
where $\tau_1$ is the unique tracial state on ${\cal  U}$ (by
\ref{meas}).

It follows from \cite{LnFU} (see 4.4.1 of \cite{Lnbk})  that there is a continuous path of unitaries $\{V(t): t\in [0,1]\}$ in
$1\otimes {\cal  U}$ such that
\beq\label{lh-5+}
V(0)=1,\,\,\, V(1)=1\otimes u_0\andeqn \text{Length}(\{V(t)\})\le
\pi +\ep/4.
\eneq

Let $\ep>0$ and ${\cal  F}\subset C\rtimes_\af \Z$ be a finite
subset.  Without loss of generality, we may assume that ${\cal
F}={\cal F}_0\cup\{u_\af\}$ for some finite subset ${\cal
F}_0\subset C.$ Define $C_1=C\otimes C(\T).$ Define $\af'\in
Aut(C_1)$ by $\af'(c\otimes 1)=\af(c)\otimes 1$ and $\af'(1\otimes
f)=1\otimes f$ for all $f\in C(\T).$ In other words,
$C_1\rtimes_{\af'} \Z \cong (C\rtimes_\af\Z)\otimes C(\T).$
Define ${\cal
F}'={\cal  F}_0\cup\{z\}.$

Define $h: (C\rtimes_\af\Z)\otimes C(\T)\to A\otimes {\cal U}$ by $h(a)=\phi(a)\otimes 1_{\cal U}$ for all $a\in C\rtimes_\af\Z$ and
$h(1\otimes z)=1\otimes u_0.$
Put ${\bar v}=v\otimes u_0.$

Let $\dt_1>0$ (in place of $\dt$), $\eta_1>0$ (in place of $\eta$), ${\cal G}_1\subset C\rtimes_\af\Z$ (in place of ${\cal G}$) be a finite subset, ${\bar {\cal G}}\subset C\otimes C(\T)$ be a finite subset and
${\cal P}\subset K_1(C\rtimes_\af\Z)$ be a finite subset required by \ref{S1uni} for $\sin(\ep/4)$ and ${\cal F}$ above.
Without loss of generality, we may assume that ${\cal G}_1={\cal G}_1'\cup \{u_\af\}$ for some finite subset
${\cal G}_1'\subset C,$ and  we may assume that ${\bar {\cal G}}={\cal G}_2'\otimes {\cal G}_2'',$ where
$1_C\in {\cal G}_2'\subset C$ is a finite subset and $\{1_{C(\T)}, z\}\subset {\cal G}_2''\subset C(\T)$ is a finite subset.

By \ref{APP}, there is $0<\dt<\dt_1$ and
 a finite subset ${\cal  G}_3\supset {\cal  G}_1'$ so that there exits a unital
\morp\, $L: (C\rtimes_{\af}\Z)\otimes C(\T) \to A\otimes {\cal  U}$ such that
\beq\label{nlh-2-1}
L(c\otimes f)\approx_{\dt_1/2} \phi(c)f(v')\rforal c\in {\cal G}_1\cup {\cal G}_2'\andeqn f\in {\cal G}_2.
\eneq
for any unitary $v'$ so that $\|[\phi(a), \, v']\|<\dt$ for all $a\in {\cal  G}_3.$
Without loss of generality, we may also assume that
\beq\label{nlh-2-1+}
[L]|_{\cal  Q}=[L']|_{\cal  Q}
\eneq
provided that $L, L'$ are both $\dt_1$-${\cal  G}_1$-multiplicative and
$L\approx_{\dt_1} L'$ on ${\cal  G}_1.$


  Let $\tau_0$ be the tracial state on
$C_1$  defined by
\beq\label{lh-2+}
\tau_0(c\otimes f)=\tau(\phi(c))\tau_1(f(u_0))\rforal c\in C\andeqn
f\in C(\T).
\eneq

Let ${\cal G}={\cal G}_1\cup {\cal G}_2'.$ Suppose that (\ref{Lh1}) holds for the above $\dt,$ ${\cal G}$ and ${\cal P}.$

Then, we have that
\beq\label{nhl-10}
\|[h(a),\, {\bar v}]\|<\dt_1\rforal a\in {\cal G}_1, \,\,\, \text{bott}_1(h|_{C\rtimes_\af\Z},\, {\bar v})|_{\cal P}=0.
\eneq
Define $L_0: C\otimes C(\T)\to A\otimes{\cal U}$ by $L_0(c\otimes f)=h(c)f({\bar v})$ for all $c\in C$ and $f\in C(\T).$
It follows from \ref{2meas} that
\beq\label{nhl-11}
\tau\circ L_0=\tau_0,
\eneq
where $\tau$ is regarded as the tracial state on $A\cong A\otimes {\cal U}.$
There exists a \morp\, $L: C\otimes C(\T)\to A\otimes {\cal U}$ such that (\ref{nlh-2-1}) holds
In particular,
\beq\label{nhl-12}
\tau\circ L\approx_{\eta} \tau\circ h|_{C_1}\,\,\,\text{on}\,\,\, {\bar {\cal G}}.
\eneq
It follows from \ref{S1uni} that there exists a unitary $W\in U(A\otimes {\cal U}\otimes {\cal U})$ such that
\beq\label{nhl-13}
\|[\psi(a),\, W]\|<\sin(\ep/4)\rforal a\in {\cal F}\andeqn W^*({\bar v}\otimes 1_{\cal U})W\approx_{\sin(\ep/4)}h(1\otimes z)\otimes 1_{\cal U},
\eneq
where $\psi: C\rtimes_\af\Z\to A\otimes {\cal U}\otimes {\cal U}$
defined by $\psi(c)=h(c)\otimes 1_{\cal U}$ for all $c\in C$ and
$\psi(u_\af)=h(u_\af)\otimes \imath(u_\af).$

In particular,
\beq\label{nlh-14}
W(1\otimes u_0\otimes 1_{\cal U})W^*\approx_{\sin(\ep/4)} v\otimes u_0\otimes 1_{\cal U}.
\eneq

By (\ref{nlh-14}),
there exists a
continuous path $\lambda: [1/4, 1/2]\to U(A\otimes {\cal U} \otimes
{\cal U} )$ such that
\beq\label{nlh-23}
\lambda(1/4)&=&v\otimes u_0\otimes 1_{\cal U},\\
\lambda(1/2)&=& W(1\otimes u_0\otimes 1_{\cal U})W^*\\
\andeqn\, \text{Length}(\{\lambda(t)\})&<&\ep/4
\eneq
Moreover,
\beq\label{lh-24}
\|[\psi(a), \lambda(t)]\|<\dt+{\ep\over{4}}<\ep\rforal a\in {\cal
F} \andeqn \rforal t\in [1/4, 1/2].
\eneq
Now define
\beq\nonumber
v_t=\begin{cases} (v\otimes 1_{\cal U})V(4t)\otimes 1_{\cal  U }, & \text{if $t\in [0,1/4)$;}\\
                              \lambda(t), &\text{if $t\in [1/4, 1/2)$ and}\\

 W( (V(2-2t))\otimes 1_{\cal U})W^*
 ,
 &\text{if $t\in [1/2,1]$.}
                             \end{cases}
                             \eneq
Then
\beq\label{lh-24+}
v_0=v\otimes 1_{\cal U} \otimes 1_{\cal  U },\,\,\, v_1=1\andeqn
\text{Length}(\{v_t\})<\pi+{\ep\over{4}}+{\ep\over{4}}+\pi+{\ep\over{4}}<2\pi+\ep.
\eneq
 For
$t\in [0, 1/4),$ we estimate that
\beq\label{lh-25}
\|[\psi(a), \, v_t]\|<\dt<\ep\rforal a\in {\cal  F}.
\eneq

By (\ref{nhl-13}), for all $c\in {\cal  F}_0$ and $ t\in [1/2,1],$
\beq\nonumber
\|[\psi(c),\, v_t]\|
&=&\|[
W^*(\phi(c)\otimes 1_{\cal U} \otimes 1_{
\cal U} )W,\, V(2-2t)\otimes 1_{\cal  U}]\|\\\label{lh-26}
&\approx_{\ep/4}&\|[ \phi(c)\otimes 1_{\cal U}\otimes 1_{\cal U} , \,V(2-2t)\otimes 1_{\cal U}]\|=0.
\eneq
 Note
again
\beq\label{lh-27}
\|[\psi(u_\af),\,V(2(t-1/2)\otimes 1_{\cal U}]\|=0.
\eneq

 Combining (\ref{nhl-13}) and (\ref{lh-27}), for $t\in [1/2, 1],$
\beq\label{lh-28}
\hspace{-0.6in}\|[\psi(u_\af),\, v_t]\| &=& \|[
W^*(\phi(u_\af)\otimes 1_{\cal U}\otimes
\imath(u_\sigma))W,\,V(2-2t)\otimes 1_{\cal U}]\|\\
&\approx_{\ep/4} &\|[\phi(u_\af)\otimes 1_{\cal U}\otimes
\imath(u_\sigma), \,V(2-2t)\otimes 1_{\cal U }]\|=0.
\eneq
It follows that
\beq\label{lh-29}
\|[\psi(c),\, v_t]\|<\ep\rforal c\in {\cal  F}\andeqn \tforal t\in
[0,1].
\eneq

\end{proof}

\section{Asymptotic unitary equivalence}

\begin{thm}\label{TAsy}
Let $C$ be a unital AH-algebra
 and let $\af\in Aut(C)$ be an
automorphism. Let $A\cong A\otimes {\cal U}$ be a unital simple
AF-algebra with a unique tracial state $\tau$ and with
$K_0(A)=\rho_A(K_0(A)).$ Suppose that $\phi_1, \phi_2:
C\rtimes_{\af}\Z\to A$ are two unital monomorphisms such that
\beq\label{Asy1+1}
\tau\circ \phi_1=\tau\circ \phi_2\andeqn {\tilde \eta}(\phi_1,
\phi_2)&=&0.
\eneq

Then, there exists a continuous path of unitaries $\{w_t: t\in [0,
\infty)\}\subset U(B\otimes {\cal  U} \otimes {\cal  U} )$ such that
\beq\label{Asy2}
\lim_{t\to\infty}{\rm ad}\, w_t\circ \psi_1''(a)=\psi_2''(a)\tforal
a\in C\rtimes_\af\Z,
\eneq
where $B=B(\psi_2, A\otimes {\cal U})$  and where $\psi_i:
C\rtimes_\af\Z\to A\otimes {\cal U} $ is defined by
$\psi_i(c)=\phi_i(c)\otimes 1_{{\cal U}}$ for all $a\in C$ and
$\psi_i(u_\af)=\phi_i(u_\af)\otimes \imath(u_\sigma),$ and
$\psi_i'': C\rtimes_\af\Z\to B\otimes {\cal U}\otimes{\cal U} $ is
defined by $\psi_i''(c)=\phi(c)\otimes 1_{{\cal  U} }\otimes
1_{{\cal  U}}$ for $c\in C$ and
$\psi_i''(u_\af)=\psi(u_\af)\otimes
\imath(u_\af),$ $i=1,2.$
\end{thm}

\begin{proof}
Since ${\rm ker}\rho_A=\{0\},$ (\ref{Asy1+1}) implies that
\beq\label{asy-1-1}
(\phi_1)_{*0}=(\phi_2)_{*0}.
\eneq
Since $K_1(A)=\{0\}$ and $K_0(A)$ is torsion free and divisible,
(\ref{asy-1-1}) implies that
\beq\label{asy-1-}
[\phi_1]=[\phi_2]\,\,\,\text{in}\,\,\, KK(C\rtimes_\af\Z, A).
\eneq

Therefore, since ${\tilde \eta}_{\phi_1, \phi_2}=0,$  there exists $\theta: K_1(C\rtimes_\af\Z) \to
K_1(M_{\phi_1, \phi_2})$ such that
\beq\label{asy-1}
(\pi_0)_{*1}\circ \theta=({\rm id}_{C\rtimes_\af \Z})_{*1}.
\eneq
Moreover,
\beq\label{asy-2}
R_{\phi_1, \phi_2}(\theta(K_1(C\rtimes_\af\Z)))=0.
\eneq
Denote by $j_0: C\to C\rtimes_\af\Z$ the natural embedding. Then one has
\beq\label{asy-3}
R_{\phi_1,\phi_2}(\theta\circ (j_0)_{*1}(K_1(C)))=0.
\eneq
Therefore
\beq\label{asy-3+}
{\tilde \eta}_{\phi_1\circ j_0, \phi_2\circ j_0}=0.
\eneq
It follows from \ref{2Tuni} that there exists a sequence of
unitaries $\{v_n\}\subset A\otimes {\cal U} $ such that
\beq\label{asy-4}
\lim_{n\to\infty}{\rm ad}\, v_n\circ\psi_1(a)=\psi_2(a)\rforal
a\in C\rtimes_\af\Z,
\eneq
where $\psi_i: C\rtimes_\af\Z\to A\otimes {\cal  U}$ defined by
$\psi_i(c)=\phi_i(c)\otimes 1_{{\cal  U}}$ for all $c\in C$ and
$\psi_i(u_\af)=\phi_i(u_\af)\otimes \imath(u_\af).$

Let $B=B(\psi_2,A\otimes {\cal  U} )$ be as in \ref{prebot}.

 To
simplify notation, without loss of generality, in what follows we
will use $\tau$ for the unique tracial state on $A,$ $A\otimes
{\cal  U}$ ($\cong A$), $B\otimes{\cal  U}$ and $B\otimes{\cal U}
\otimes {\cal U} \cong B\otimes {\cal  U} .$

Let $\{{\cal  F}_n\}\subset C\rtimes_\af\Z$ be an increasing
sequence of finite subsets of the unit ball of $C\rtimes_\af\Z$
whose union is dense in the unit ball. Let $\dt_n>0,$ ${\cal
G}_n\subset C\rtimes_{\af}\Z$ be a finite subset and ${\cal
P}_n\subset K_1(C\rtimes_\af\Z)$ be a finite subset required by
\ref{Lhom} for $\ep_n=1/2^{n+2}$ and ${\cal  F}_n.$ We may assume
that $\text{bott}_1(h', v')|_{{\cal  P}_n}$ is well defined for any
unital \hm\, $h'$ and any unitary $v',$  whenever
\beq\label{asy-5}
\|[h'(c),\, v']\|<\dt_n\rforal c\in {\cal  G}_n
\eneq
for $n=1,2,....$
We also assume that
\beq\label{asy-5+1}
\text{bott}_1(h',v')|_{{\cal  P}_n}=\text{bott}_1(h', v'')|_{{\cal  P}_n}
\eneq
provided that $\|v'-v''\|<\dt_n.$

We may assume that $\dt_n<1/2^{n+1}$ and ${\cal
F}_n\subset {\cal  G}_n,$  $n=1,2,....$
We may write ${\cal  P}_n=\{x_1,x_2,...,x_{r(n)}\},$ $n=1,2,....$

Let $\dt_n'>0$ and $\{y_1,y_2,...,y_{m(n)}\}\subset
K_1(C\rtimes_\af\Z)$ be a finite subset required by \ref{Bott1}
for $\dt_n/2,$ ${\cal  G}_n$ and ${\cal  P}_n.$ Put ${\cal  Q}_n=\{y_1,y_2,...,y_{m(n)}\}.$
We may assume that $m(n+1)\ge m(n).$
To simplify notation, we may assume that $x_j=y_j,$ $j=1,2,...,r(n)$ and $r(n)\le m(n),$
$n=1,2,....$

There is an integer $l(n)\ge 1$ and unitaries
$z_1,z_2,...,z_{m(n)}\in M_{l(n)}(C\rtimes_\af\Z)$ such that
$[z_j]=y_j,$ $j=1,2,...,m(n).$

Let $\eta_n>0$ be such that $\eta_n<\min\{\dt_n/2l(n)^2, (1/4)\sin(\dt_n'/2l(n)^2)\}$ and
$\eta_{n+1}<\eta_n/2,$
$n=1,2,....$
We may also assume that, by passing to a subsequence, if necessary
\beq\label{asy-6}
{\rm ad}\, v_n\circ \psi_1\approx_{\eta_n/2}\psi_2\,\,\,\text{on}\,\,\, {\cal  G}_n
\eneq
Note that
\beq\label{asy-7}
\|[\psi_1(a),\, v_nv_{n+1}^*]\|<\eta_n\rforal a\in {\cal  G}_n
\eneq
So $\text{bott}_1(\phi_1,v_nv_{n+1}^*)|_{{\cal  Q}_n}$ is well-defined.

Put
$$
{\bar v}_n={\rm diag}(\overbrace{v_n, v_n,...,v_n}^{l(n)}).
$$
 We will continue to use $\phi_i$ for $\phi_i\otimes_{{\rm id}_{M_{l(n)}}}$ and use $\psi_i$ for
 $\psi_i\otimes_{{\rm id}_{M_{l(n)}}},$ $i=1,2,$ respectively.
 Thus, we may assume that, for $j\le m(n),$
 \beq\label{asy-8}
&& \|\psi_2(z_j){\rm ad}\, {\bar v}_n\circ \psi_1(z_j^*)-1\|<(1/4)\sin (\dt_n'/2)
 \eneq
Put
\beq\label{asy-9}
h_{j,n}=\log({1\over{2\pi i}}\psi_2(z_j){\rm ad}\, {\bar v}_n\circ \psi_1(z_j^*)),\,\,\,j=1,2,...,m(n),\,n=1,2,....
\eneq
Then
\beq\label{asy-10}
\tau(h_{j,n})<\dt_n/2,
\eneq
By (\ref{asy-3}), $R_{\psi_1,\psi_2}(K_1(M_{\psi_1,
\psi_2}))\subset \rho_{A\otimes {\cal U}}(K_0(A\otimes {\cal U} )).$ It follows from 3.5 of \cite{Lnasy} that
\beq\label{asy-11}
\widehat{h_{j,n}}(\tau)=\tau(h_{j,n})\in \rho_{A\otimes {\cal U}}(K_0(A\otimes {\cal  U})),\,\,\,j=1,2,...,m(n), n=1,2,....
\eneq
Since $K_0(A\otimes {\cal U} )$ is divisible, there exists a \hm\, $\bt_n: K_1(C\rtimes_\af\Z)\to K_0(A\otimes {\cal  U})$ (see also 6.1, 6.2 and 6.3 of \cite{Lnemb2})
such that
\beq\label{asy-12}
\tau(\bt_n(y_j))={\widehat{h_{j,n}}}(\tau)=\tau(h_{j,n}),
j=1,2,...,m(n), n=1,2,....
\eneq
It follows from \ref{Bott1} that there exists a unitary $U_n\in
B=B(\psi_2, A\otimes {\cal U} )$ such that
\beq\label{asy-23}
\|[\psi_2'(a),\, U_n]\|<\dt_n'\rforal a\in {\cal  G}_n\andeqn
\text{bott}_1(\psi'_2,\, U_n)|_{{\cal  P}_n}=-\bt_n'|_{{\cal
P}_n},
\eneq
where $\bt_n'=(j_2)_{*0}\circ \bt_n,$ $n=1,2,...,$ and $j_2:
A\otimes {\cal  U} \to B$ is the embedding and where $\psi_i:
C\rtimes_\af\Z\to B\otimes {\cal  U}$ is defined by
$\psi_i'(c)=\psi(c)\otimes 1_{{\cal U} }$ for all $c\in
C\rtimes_\af\Z$ and $\psi(u_\af)\otimes \imath(u_\af).$ Note that
(see \ref{prebot}) that $\tau\circ j_2(a)=\tau(a)$ for $a\in
A\otimes {\cal U}.$

By Exel's trace formula (see Theorem 3.6 of \cite{Lnasy}) for $j=1,2,...,r(n),$
\beq\label{asy-24}
\tau(h_{j,n})&=&-\text{bott}_1(\psi_2',\, U_n)(x_j)(\tau)\\
  &=& -\tau(\log({1\over{2\pi i}}{\bar U}_n\psi_2'(z_j){\bar U}_n^*\psi_2'(z_j^*))),
  \eneq
where
$$
{\bar U}_n={\rm diag}(\overbrace{U_n,U_n,...,U_n}^{l(n)}).
$$
Define
$u_n=v_nU_n,$ $n=1,2,....$ Put
\beq\label{asy-25}
{\bar u}_n={\rm diag}(\overbrace{u_n,u_n,...,u_n}^{l(n)})\andeqn
{\bar u}_{n+1}'={\rm diag}(\overbrace{u_{n+1},u_{n+1},...,u_{n+1}}^{l(n)}),
\eneq
$n=1,2,....$
By 6.1 of \cite{Lnemb2} and (\ref{asy-24}), we compute that
\beq\label{asy-26}
&&\hspace{-0.3in}\tau(\log({1\over{2\pi i}}(\psi'_2(z_j){\rm ad}\, {\bar
u_n}(\psi'_1(z_j^*)))))\\
&=&\tau(\log({1\over{2\pi i}}({\bar U_n}\psi_2'(z_j^*){\bar
U_n}^*{\bar v}_n^*
(\psi'_1(z_j)){\bar v_n})))\\
&=&\tau(\log({1\over{2\pi i}}({\bar U_n}\psi'_2(z_j){\bar
U_n}^*\psi'_2(z_j^*)\psi_2'(z_j){\bar v}_n^*
(\psi'_1(z_j^*)){\bar v_n})))\\
&=&\tau(\log({1\over{2\pi i}}({\bar U_n}\psi'_2(z_j){\bar
U_n}^*\psi'_2(w_j^*))))\\
&&\hspace{0.6in}+\tau({1\over{2\pi i}}\log(\psi'_2(z_j){\bar v}_n^*
(\psi'_1(z_j^*)){\bar v_n})))\\\label{asy-26+}
&=&\text{bott}_1(\psi'_2, U_n))(z_j)(\tau)+\tau(h_{j,n})=0.
\eneq

Let
\beq\label{asy-27}
b_{j,n}&=&\log({1\over{2\pi i}}{\bar u_n}\psi'_2(z_j){\bar
u_n}^*\psi_1'(z_j^*))\andeqn\\
 b_{j,n}'&=&\log({1\over{2\pi
i}}\psi'_2(z_j){\bar u_n}^*{\bar u_{n+1}}'\psi'_2(z_j^*)({\bar
u_{n+1}}')^*{\bar u_n})),
\eneq
$j=1,2,...,n$ and $n=1,2,....$
 We have, by (\ref{asy-26+}),
\beq\label{asy-28}
\tau(b_{j,n})&=&\tau(\log({1\over{2\pi i}}{\bar
u_n}\psi'_2(z_j){\bar
u_n}^*\psi_1'(z_j^*)))\\
&=& \tau(\log{1\over{2\pi i}}{\bar u_n}^*{\bar
u_n}\psi'_2(z_j){\bar u_n}^*\psi'_1(z_j^*){\bar u_n})\\
&=&\tau(\log{1\over{2\pi i}}\psi'_2(z_j){\bar u_n}^*\psi'_1(z_j^*)
{\bar u_n})=0.
\eneq
 $j=1,2,...,r(n)$ and $n=1,2,....$  Note also
$\tau(b_{j,n+1})=0$ for $j=1,2,...,r(n).$ In $M_{l(n+1)}(C\rtimes_{\af}\Z),$ $x_j$
represented by
\beq\label{asy-28+}
{\rm diag}(z_j, 1_{M_{l(n+1)-l(n)}})
\eneq
for $j=1,2,...,r(n).$
Put
\beq\label{asy-28++}
b_{j,n+1}''=\log({1\over{2\pi i}}\psi'_2(z_j)({\bar
u_{n+1}}')^*\psi'_1(z_j^*){\bar u_{n+1}}')
\eneq
Therefore
\beq\label{asy-29-1}
\tau(b_{j,n+1}'')=\tau(b_{j, n+1})=0.
\eneq
 We have that
\beq\label{asy-29}
{\bar u_n}^*e^{2\pi i b_{j,n}'}{\bar u_n}=e^{2\pi i
b_{j,n}}e^{-2\pi i b_{j,n+1}''}.
\eneq
Thus, by 6.1 of \cite{Lnemb2}, we compute that
\beq\label{asy-30}
\tau(b_{j,n}')=\tau(b_{j,n})-\tau(b_{j,n+1})=0.
\eneq
It follows the Exel trace formula (Theorem 3.6 of \cite{Lnasy}) and (\ref{asy-30})
that
\beq\label{asy-31}
&&\hspace{-0.3in}\text{bott}_1(\psi'_2,u_n^*u_{n+1}))(x_j)(\tau)\\
&=&-\tau(\log({\bar u_n}^*{\bar
u_{n+1}}'\psi_2(z_j^*)({\bar u_{n+1}}')^*{\bar u_n}\psi_2(z_j))\\
&=&-\tau(\log(\psi'_2(z_j){\bar u_n}^*{\bar
u_{n+1}}'\psi'_2(z_j^*)({\bar u_{n+1}}')^*{\bar u_n}))=0.
\eneq

By applying \ref{Lhom}, we obtain a continuous path of unitaries $\{v_n(t): t\in [0,1]\} \subset B\otimes {\cal U}\otimes{\cal  U} $
 such that
 \beq\label{asy-32}
 v_n(0)=u_nu_{n+1}^*,\,\,\, v_n(1)=1\andeqn \|[\psi_2''(a),\,v_n(t)]\|<1/2^{n+2}\rforal a\in {\cal  F}_n,
 \eneq
where $\psi_i'': C\rtimes_\af\Z\to B\otimes{\cal U}\otimes {\cal  U} $ is defined by
$\psi_i''(c)=\psi_i'(c)$ for all $c\in C$ and $\psi_i''(u_\af)=\psi_i'(u_\af)\otimes\imath(u_\sigma),$ $i=1,2$ and
 $n=1,2,....$
 Define $u(t)=v_n(t-n+1)u_{n+1}$ for $t\in [n-1, n),$ $n=1,2,....$
We check that
\beq\label{asy-33}
\lim_{t\to\infty}{\rm ad}\, u(t)\circ \psi''_1(a)=\psi''_2(a)\rforal a\in C\rtimes_\af\Z.
\eneq

\end{proof}

\section{AF-embedding of $C\rtimes_{\Lambda}\Z^2$}

\begin{thm}\label{1MT}
Let $C$ be a unital AH-algebra and let $\Lambda: \Z^2\to Aut(C)$
be a   \hm. Then $C\rtimes_{\Lambda}\Z^2$ can be embedded
into a unital simple AF-algebra if and only if $C$ admits a
faithful $\Lambda$-invariant tracial state.

\end{thm}

\begin{proof}
We only need to show the ``if " part of the theorem. Let $\af_1$
and $\af_2$ be two generators of $\Lambda(\Z^2).$ Thus $\af_1,
\af_2\in Aut(A),$ $\af_1\circ \af_2=\af_2\circ \af_1$ and
$C\rtimes_{\Lambda}\Z^2\cong(C\rtimes_{\af_1}\Z)\rtimes_{\af_2}\Z.$

Put $C_1=C\rtimes_{\af_1}\Z.$ We view $\af_2$ as an automorphism
on $C_1$ by $\af_2(u_{\af_1})=u_{\af_1}.$

Let $t$ be a faithful $\Lambda$-invariant tracial state on $C.$ We
may also view $t$ as a faithful tracial state on $C_1$ as well as
on $C\rtimes_{\Lambda}\Z^2.$

It follows from (the proof) of 10.5 of \cite{Lnasy} that there exists
a unital simple AF-algebra $A\cong A\otimes {\cal U}$ with unique tracial state $\tau$
and with $K_0(A)=\rho_A(K_0(A))$ and a unital
monomorphism $h: C_1\to A$ such that
\beq\label{1mt-1}
\tau\circ h=t.
\eneq
Put $\phi_1=h$ and $\phi_2=h\circ \af_2.$ Then
\beq\label{1mt-2}
\tau\circ \phi_1=\tau\circ \phi_2.
\eneq
In particular, since ${\rm ker}\rho_A=\{0\},$ $K_1(A)=\{0\}$  and
$K_0(A)$ is torsion free and divisible,
\beq\label{lmt-3}
[\phi_1]=[\phi_2]\,\,\,\text{in}\,\,\, KK(C_1, A).
\eneq

There exists $\theta: K_1(C_1)\to K_1(M_{\phi_1,\phi_2}).$ Let
$\D_0=R_{\phi_1,\phi_2}(\theta(K_1(C_1))).$ Let $\D_1$ be a
countable subring of $\R$ containing $\rho_A(K_0(A))\cup \D_0$
which is also a divisible subgroup of $\R.$ Let $A_1$ be a unital
simple AF-algebra with a unique tracial state such that
$K_0(A_1)=\D_1.$ Moreover, there is a unital embedding $j: A\to A_1$
such that
\beq\label{lmt-5}
\tau\circ j=\tau,
\eneq
where we also use $\tau$ for the unique tracial state on $A_1.$ To
simplify notation, without loss of generality, we may assume that
$A=A_1.$ Therefore, with this assumption, we have
\beq\label{lmt-6}
R_{\phi_1,\phi_2}(M_{\phi_1, \phi_2})\subset \rho_A(K_0(A))=\D_1.
\eneq
By 10.3 of \cite{Lnasy}, this implies that
\beq\label{lmt-7}
{\tilde \eta}_{\phi_1,\phi_2}=0.
\eneq

 It follows that \ref{TAsy} that there exists a continuous path of
 unitaries $\{w_t: t\in [0, \infty)\}$ of $B\otimes
 {\cal U} \otimes {\cal U} )$ such that
 \beq\label{lmt-7+}
 \lim_{t\to\infty}{\rm ad}\, w_t\circ
 \psi_1''(a)=\psi_2''(a)\tforal a\in C_1,
 \eneq
where $\psi_i'': C_1\to B\otimes {\cal U} \otimes {\cal  U} $ is
defined by $\psi_i''(c)=\phi(c)\otimes 1_{{\cal  U}}\otimes 1_{{\cal
U} }$ for all $c\in C$ and
$\psi_i''(u_{\af_1})=\phi_i(u_{\af_1})\otimes
\imath(u_\sigma)\otimes \imath(u_\sigma),$ and $B=B(\psi_2, A\otimes
{\cal  U} )$ and $\psi_2: C_1\to A\otimes {\cal  U}$ is defined by
$\psi_2(c)=\phi_2(c)\otimes 1_{{\cal  U} }$ for all $c\in C$ and
$\psi_2(u_{\af_1})=\phi_2(u_{\af_1})\otimes\imath(u_{\sigma})\otimes
\imath(u_{\sigma}).$ It follows from 5.5 of \cite{Lnemb2} that there
exists a unital monomorphism $h_1: C_1\rtimes_{\af_2}\Z\to B_1,$
where $B_1$ is a unital simple AF-algebra with a unique tracial
state.

\end{proof}

\begin{cor}\label{Tz2}
Let $C$ be a unital AH-algebra and let $\Lambda: \Z^2\to Aut(C)$
be a \hm. Suppose that $t$ is a faithful
$\Lambda$-invariant tracial state on $C.$

Then there exists a unital simple AF-algebra $B_1$ with a unique
tracial state $\tau$ and with
$K_0(B_1)=\rho_{B_1}(K_0(B_1))$ which is also a subring of $\R,$
and there exists a unital monomorphism $h_1:
C\rtimes_{\Lambda}\Z^2\to B_1$ such that
\beq\label{Tz2-1}
\tau\circ h_1=t.
\eneq

\end{cor}

\begin{proof}
At the end of the proof of \ref{1MT}, we apply 5.5 of \cite{Lnemb2}.
The proof of 5.5 of \cite{Lnemb2} actually shows that one can have
the further requirements for $B_1,$ $h_1$ and $\tau.$

\end{proof}

\begin{cor}\label{CV}
Let $X$ be a compact metric space and $\Lambda: \Z^2\to Homeo(X)$
be a  \hm. Then $C(X)\rtimes_{\Lambda}\Z^2$ can be
embedded into a unital simple AF-algebra if and only if $X$ admits
a strictly positive $\Lambda$-invariant probability Borel measure.
\end{cor}

\section{The absorption lemma}

\begin{df}\label{sigmak}
{\rm

Let $C$ be a unital separable amenable \CA\, and let $\Lambda_k:
\Z^k\to Aut(C)$ be a \hm. Suppose that
$\{\af_1,\af_2,...,\af_k\}$ is the set of generators of
$\Lambda(\Z^k).$

Suppose that $\phi: C\rtimes_{\Lambda}\Z^k\to A$ is a unital
monomorphism. Denote by $\phi^{(k)}: C\rtimes_{\Lambda_1}\Z\to
A\otimes {\cal U}^k$ the monomorphism defined by
\beq
\phi^{(k)}(c)&=&\phi(c)\otimes 1_{{\cal U}^{k}}\hspace{0.8in}\rforal c\in C,\\
\phi^{(k)}(u_{\af_1})&=&\phi(u_{\af_1})\otimes
\imath(u_\sigma)
\otimes 1_{{\cal U}^{k-1}},\\
\phi^{(k)}(u_{\af_2})&=& \phi(u_{\af_2})\otimes 1_{\cal
U}\otimes \imath(u_\sigma)\otimes 1_{{\cal U}^{k-1}},\\
&\cdots &\cdots\\
\phi^{(k)}(u_{\af_{k}})&=&(\phi(u_{\af_k})\otimes 1_{{\cal
U}^{(k-1)}})\otimes \imath(u_\sigma).
\eneq

Let $\{\af_1,\af_2,...,\af_k\}$ be a set of generators for $\Lambda(\Z^k).$ Define
$\Lambda_j:\Z^j\to Aut(C)$ by $\{\af_1,\af_2,...,\af_j\},$ $j=1,2,...,k.$
Denote by $\phi_j: C\rtimes_{\Lambda_j}\Z^j\to A$ the \hm\, $\phi|_{C\rtimes_{\Lambda_j}\Z^j}$ and
$\phi^{(j)}=\phi_j^{(j)},$ $j=1,2,...,k.$

Denote by  $B_1=B(\phi, A)$ the AF-algebra defined in \ref{prebot}. Define $B_{j+1}(\phi, A)=B(\phi^{(j+1)}, B_j\otimes {\cal U}^j),$
$j=1,2,...,k-1.$

 }

\end{df}

\begin{lem}\label{rot}
In ${\cal U},$ there is a continuous and piecewise smooth path of unitaries $\{u(t):t\in [0, 1]\}$ in ${\cal U}$ such that
\beq\label{rot1}
&&u(0)=\imath(u_\sigma),\,\,\, u(1)=1_{\cal U}\andeqn\\
&&\int_0^1 \tau({du(t)\over{dt}}u(t)^*)dt=0,
\eneq
where $\tau$ is the unique tracial state of ${\cal U}.$
\end{lem}

\begin{proof}
There are two self-adjoint elements $h_1, h_2\in {\cal U}$ such that
\beq\label{rot-1}
\phi(u_\sigma)=e^{2\pi ih_1}e^{2\pi ih_2},
\eneq
where $\phi: {\cal U}\rtimes_{\sigma}\Z\to {\cal U}$ is defined in \ref{U1} (see \cite{LnFU}).

Define $v(t)=e^{2\pi ih_1(1-t)}e^{2\pi (i h_2(1-t)}.$ Then $v(0)=\phi(u_\sigma)$ and
$v(1)=1_{\cal U}.$
Define $w(t)=v(t)\otimes v(t)^*$ for $t\in [0,1].$
Then,
\beq
{dw(t)\over{dt}}w(t)^*&=& {d((v(t)\otimes 1_{\cal U})(1_{\cal U}\otimes v(t)^*)\over{dt}}((v(t)^*\otimes 1_{\cal U})
(1_{\cal U}\otimes v(t)))\\
&=&(-ih_1)\otimes 1_{\cal U}+(-ie^{2\pi ih_1(1-t)}h_2e^{-2\pi i h_1(1-t)})\otimes 1_{\cal U}\\
&+& 1_{\cal U}\otimes (ih_1)\otimes 1_{\cal U}+(ie^{2\pi ih_1(t-1)}h_2e^{2\pi ih_1(1-t)})\otimes 1_{\cal U}
\eneq
Therefore
\beq\label{rot-3}
\tau({dw(t)\over{dt}}w(t)^*)=0.
\eneq
Define $u(t)=s(w(t))$ for $t\in [0,1],$ where $s$ is defined in \ref{U1}.
Note that $u(0)=s(w(0))=s(\phi(u_\sigma)\otimes \phi(u_\sigma)^*)=\imath(u_\sigma)$ and
$u(1)=1_{\cal U}.$ Moreover,
\beq\label{rot-4}
\tau({du(t)\over{dt}}u(t)^*)&=&\tau(s({dw(t)\over{dt}}w(t)^*))\\
&=&\tau({dw(t)\over{dt}}w(t)^*)=0
\eneq
for all $t\in [0,1].$

\end{proof}

\begin{lem}\label{uuu}
There exists a continuous  path of unitaries $\{w(t):t\in [0, \infty)\}$ of ${\cal U}\otimes {\cal U}$
such that
\beq\label{uuu1}
\lim_{t\to\infty} w(t)^*(\imath(u_\sigma)\otimes \imath(u_\sigma))w(t)=\imath(u_\sigma)\otimes 1_{\cal U}.
\eneq
\end{lem}

\begin{proof}
Let $\tau$ be the unique tracial state of ${\cal U}.$ Using the standard conditional expectation, one computes that
\beq\label{uuu-1}
\tau(\imath(u_{\sigma}))=0.
\eneq
Let $f=\sum_{i=1}^n\lambda_i u_\sigma^i$ be a polynomial of $u_\sigma,$ where $\lambda_i\in \C.$
Then
\beq\label{uuu-2}
\tau(\imath(f))=\lambda_0.
\eneq
Similarly
\beq\label{uuu-3}
\tau(\sum_{i=1}^n \lambda_i(\imath(u_\sigma)\otimes\imath(u_\sigma))^i)=\lambda_0.
\eneq
It follows that
\beq\label{uuu-4}
\tau(f(\imath(u_\sigma)))=\tau(f(\imath(u_\sigma)\otimes \imath(u_\sigma)))
\eneq
for all $f\in C(\T),$ where we also use $\tau$ for the unique tracial state of ${\cal U}\otimes {\cal U}.$

Define $h_1, h_2: C(\T)\to {\cal U}\otimes {\cal U}$ by
$h_1(f)=f(\imath(u_\sigma)\otimes\imath(u_\sigma))$ and
$h_2(f)=f(\imath(u_\sigma))\otimes 1_{\cal U}$  for all $f\in C(\T).$
Since $K_1({\cal U})=\{0\},$ by (\ref{uuu-4}),
we have
\beq\label{uuu-5}
[h_1]=[h_2]\,\,\,\text{in}\,\,\,KK(C(\T), {\cal U})\andeqn
\tau\circ h_1=\tau\circ h_2.
\eneq
Define $\theta: C(\T)\to M_{h_1, h_2}$ by
\beq\label{uuu-6}
\theta(f)(t)=\imath(u_\sigma)\otimes u(t)\tforal t\in [0,1],
\eneq
where $u(t)$ is a continuous piecewise smooth path of unitaries in ${\cal U}$ such that
\beq\label{uuu-7}
u(0)=\imath(u_\sigma),\,\,\, u(1)=1_{\cal U}\andeqn
\int_0^1 \tau({du(t)\over{dt}}u(t)^*)dt=0,
\eneq
given by \ref{rot}.
Then
${\tilde \eta}_{h_1,h_2}=0.$
By applying the main theorem of \cite{Lnasy}, there exists a continuous path of unitaries
$\{w(t): t\in [0, \infty)\}\subset {\cal U}\otimes {\cal U}$ such that
\beq\label{uuu-9}
\lim_{t\to\infty}w(t)^*(\imath(u_\sigma)\otimes \imath(u_\sigma))w(t)=\imath(u_\sigma).
\eneq

\end{proof}

The proof of the following is similar to  but  easier than that of \ref{uuu}.

\begin{lem}\label{2uuu}
There exists a continuous path of unitaries $\{u(t): t\in [0, \infty)\} \subset {\cal U}\otimes {\cal U}$
such that
\beq\label{2uuu-1}
\lim_{t\to\infty}u(t)^* (\imath(u_\af)\otimes 1_{\cal U})u_t=1_{\cal U}\otimes \imath(u_\af).
\eneq
\end{lem}

\begin{lem}\label{Abs}
Let $C$ be a unital separable amenable \CA, let $\Lambda: \Z^k\to Aut(A)$ be a \hm\,
such that $\{\af_1,\af_2,...,\af_k\}$ forms a set of generators for $\Lambda(\Z^k).$
Let $\phi: C\rtimes_{\Lambda}\Z^k\to A$ be a monomorphism.
Then there exists a continuous path of unitaries $\{w(t): t\in [0, \infty)\}$ of $A\otimes {\cal U}^{2k}$
such that
\beq\label{Abs1}
\lim_{t\to\infty}{\rm ad}\, w(t)\circ (\phi^{(k)})^{(k)}(a)={\phi}^{(k)}(a)\otimes 1_{{\cal U}^k}\rforal a\in C.
\eneq

\end{lem}

\begin{proof}
Define $\psi: C\rtimes_{\Lambda}\Z^k\to A\otimes {\cal U}^{2k}$ as follows:
\vspace{-0.1in}
\beq\label{abs-1}
\psi(c)&=& \phi(c)\otimes \overbrace{1_{{\cal U}}\oplus 1_{{\cal U}}\oplus \cdots\oplus 1_{{\cal U}}}^{2k}\,\,\, \hspace{0.4in}\rforal c\in C,\\
\psi(u_{\af_1})&=&\phi(u_{\af_1})\otimes \imath(u_\sigma)\otimes \imath(u_\sigma)\otimes 1_{{\cal U}^{2k-2}},\\
\psi(u_{\af_2}) &=&\phi(u_{\af_2})\otimes 1_{{\cal U}^2}\oplus \imath(u_\sigma)\otimes \imath(u_\sigma)\otimes 1_{{\cal U}^{2k-4}},\\
&\cdots&\,\cdots\\
\psi(u_{\af_k})&=&\phi(u_{\af_k})\oplus 1_{{\cal U}^{2k-2}}\oplus \imath(u_\sigma)\otimes \imath(u_\sigma).
\eneq
It is clear, by \ref{2uuu}, that there exists a continuous path of unitaries
$\{u(t): t\in [0,\infty)\}\subset {\cal U}^{2k}$ such that
\beq\label{abs-2-1}
\lim_{t\to\infty} (1\otimes u(t))^* (\phi^{(k)})^{(k)}(a) (1\otimes u(t))=\psi(a)\rforal a\in C.
\eneq

 Therefore, to complete the proof,
it suffices to show that there exists a continuous path of unitaries $\{w(t): t\in [0, \infty)\}$ of $A\otimes {\cal U}^{2k}$
such that
\beq\label{abs-2}
\lim_{t\to\infty}{\rm ad}\, w(t)\circ \psi(a)={\phi}^{(k)}(a)\otimes 1_{{\cal U}^k}\rforal a\in C.
\eneq

But this clearly follows from \ref{uuu}.

\end{proof}

\begin{cor}\label{casy}
Let $C$ be a unital AH-algebra
 and let $\af\in Aut(C)$ be an
automorphism. Let $A\cong A\otimes {\cal U}$ be a unital simple AF-algebra with a unique
tracial state $\tau$ and with $K_0(A)=\rho_A(K_0(A)).$
 Suppose that $\phi_1, \phi_2: C\rtimes_{\af}\Z\to A$
are two unital monomorphisms such that
\beq\label{nAsy1+1}
\tau\circ \phi_1=\tau\circ \phi_2\andeqn {\tilde \eta}(\phi_1,
\phi_2)=0.
\eneq

Then, there exists a continuous path of unitaries $\{w_t: t\in [0,
\infty)\}\subset U(B\otimes {\cal  U})$ such that
\beq\label{nAsy2}
\lim_{t\to\infty}w_t^*(\psi_1(a)\otimes 1_{\cal U})w_t=\psi_2(a)\otimes 1_{\cal U}\tforal a\in
C\rtimes_\af\Z,
\eneq
where $B=B(\phi_2^{(1)}, A\otimes {\cal U})$  and with $1_A=1_B,$
and $\psi_i=j_2\circ \phi_i^{(1)}$ {\rm (}where $j_2: A\otimes
{\cal U}\to B$ is defined in {\rm \ref{prebot}}{\rm )}, $i=1,2.$

\end{cor}

\begin{proof}
This follows immediately from \ref{TAsy} and \ref{Abs}.
\end{proof}

\begin{lem}\label{3abs}
Let $B_1\cong B_1\otimes {\cal U}$ and $B_2\cong B_2\otimes {\cal U}$ be two
unital simple AF-algebras with unique tracial state and with
$\rho_{B_i}(K_0(B_i))=K_0(B_i)$ which is also a subring of $\R$
($i=1,2$). Suppose also that $B_1\otimes {\cal U}^k\subset B_2,$
 $1_{B_1\otimes {\cal U}^k}=1_{B_2}$ for some $k\ge 1$ and suppose that $j$ is the embedding
 from $B_1\otimes {\cal U}^k$ to $B_2.$  Let
$j_0: B_1\to B_2$ be a unital embedding such that $(j_0)_{*0}={\rm
j}_{*0}.$
 Then there is a continuous path of unitaries $\{u_t:
t\in [0,\infty)\}\subset B_2\otimes {\cal U}^k$ such that
\beq\label{3abs1}
\lim_{t\to\infty} u_t^* (a\otimes b)u_t=j_0(a)\otimes b
\eneq
for all $a\in B_1$ and $b\in {\cal U}^k.$

\end{lem}
\begin{proof}
Note that $[j_0]=[j]$ in $KK(B_1, B_2).$ Let $\tau$ be the unique
tracial state of $B_2$ then
$$
\tau\circ j=\tau\circ j_0.
$$
Since $K_1(B_1)=\{0\},$ ${\tilde \eta}_{j_0, j}=0.$ It follows theorem 9.1
of \cite{Lnasy} (see also Theorem 1 of \cite{M1}) that (\ref{3abs1}) holds.

\end{proof}

\section{The general cases}

In this section, we will show the main theorem of this paper.
We will repeat the same argument deployed for the proof of \ref{1MT}
and  use induction to present the proof.

\begin{df}\label{DP1}

{\rm Let $C$ be a unital separable amenable \CA\, in ${\cal  N}.$

(1)\, We say $C$ has the property (P1) with integer $k\ge 1,$ if
the following statement holds:

Let $A\cong A\otimes {\cal U}$ be a unital simple AF-algebra with a unique tracial state
$\tau$ and with $K_0(A)=\rho_A(K_0(A)).$ There exists at least one unital monomorphism
from $C\rtimes_{\Lambda}\Z^k$ into $A.$

Suppose
that $\phi_1, \phi_2: C\rtimes_{\Lambda}\Z^k\to A$ are two unital
monomorphisms such that
\beq\label{dAsy1+1}
\tau\circ \phi_1=\tau\circ \phi_2\andeqn {\tilde \eta}(\phi_1,
\phi_2)&=&0.
\eneq

Then, there exists a continuous path of unitaries $\{w_t: t\in [0,
\infty)\}\subset U(B \otimes {\cal  U}^{k}
)$ such that
\beq\label{dAsy2}
\lim_{t\to\infty}w_t^*(\psi_1(a)\otimes 1_{{\cal
U}^k})w_t=\psi_2(a)\otimes 1_{{\cal U}^k}\tforal a\in
C\rtimes_{\Lambda}\Z^k,
\eneq
where $A\subset B\cong B\otimes {\cal U}$
 is a unital simple AF-algebra with a unique
tracial state, with $K_0(B)=\rho_B(K_0(B)$ which is also a subring
of $\R$ and with $1_A=1_B,$ and where $\psi_i=j\circ \phi_i^{(k)}$
{\rm (}where $j: A\otimes {\cal U}^k\to B$ is a unital
monomorphism {\rm )}, $i=1,2.$
defined in {\rm \ref{prebot}}{\rm )} $i=1,2,$

\vspace{0.1in}

(2) We say $C$ has the property (P2) with integer $k,$ if the
following statement holds:

Let $\Lambda:\Z^k\to Aut(C)$ be a \hm, let $A\cong A\otimes {\cal
U} $ be a unital simple AF-algebra with a unique tracial state and
with $K_0(A)=\rho_A(K_0(A))$ which is also a subring of $\R.$
Suppose that there is at least one unital monomorphism from
$C\rtimes_{\Lambda}\Z^k$ to $A.$
Suppose that  $j: C\rtimes_{\Lambda}\Z^k\to  A$ is a unital
monomorphism.
Let $\ep>0,$ ${\cal  F}\subset C\rtimes_{\Lambda}\Z^k$ be a finite subset
and let ${\cal P}\subset K_1(C\rtimes_{\Lambda}\Z^k)$ be a finite subset. There is
$\dt>0$ and a finite subset $\{y_1,y_2,...,y_m\}\subset
K_1(C\rtimes_{\Lambda}\Z^k)$ satisfying the following:

 If $\bt: K_1(C\rtimes_{\Lambda}\Z^k)\to K_0(A)$ is a \hm\, with
 \beq\label{Bott1-1-}
 \rho_A((\bt(y_j))(\tau)<\dt,\,\,\,j=1,2,...,m,
 \eneq
 where $\tau$ is the unique tracial state of $A,$
 then there exists a    unitary $u\in B$ such that
 \beq\label{Bott1-2-}
 \|[\psi(a), \, u]\|<\ep\tforal a\in {\cal  F}\andeqn
 {\rm{bott}}_1(\psi,\, u)|_{{\cal  P}}=\bt'|_{{\cal  P}}
 \eneq
where $\psi=j^{(k)}$
and $\bt'=(j_2)_{*0}\circ \bt$ ($j_2: A$ is as in \ref{prebot}),
where $A\subset B\cong B\otimes{\cal U}$ is a unital simple
AF-algebra with a unique tracial state, with
$K_0(B)=\rho_B(K_0(B)$ which is also a subring of $\R$ and with
$1_A=1_B$ ($j_2: A\to B$ is the unital embedding).

\vspace{0.1in}

(3) We say that $C$ has the property (P3) with integer $k,$ if the
following statement holds:

Let $\Lambda: \Z^k\to Aut(C)$ be a \hm\, and let $A\cong A\otimes
{\cal U}$ be a unital simple AF-algebra with a unique tracial
state $\tau$ and $\rho_A(K_0(A))=K_0(A).$ There exists at least
one unital monomorphism from $C\rtimes_{\Lambda}\Z^k$ into $A.$

Let $h:
(C\rtimes_{\Lambda}\Z^k)\otimes C(\T)\to A$ be a unital
monomorphism. Let $\{\af_1,\af_2,...,\af_k\}$ be the set of
generators of $\Lambda(\Z^k)$ and defined $\Lambda_{k-1}:
\Z^{k-1}\to Aut(A)$ by $\{\af_1,\af_2,...,\af_{k-1}\}.$
For any $\ep>0$ and a finite subset ${\cal  F}\subset C\rtimes_{\Lambda}
\Z^k,$ there exists $\dt>0,$  $\eta>0,$ a finite subset ${\cal
G}\subset C\rtimes_{\Lambda}\Z^k,$  a finite subset ${\bar {\cal
G}}\subset C\rtimes_{\Lambda_{k-1}}\Z^{k-1}\otimes C(\T)$ and a
finite subset ${\cal P}\subset K_1(C\rtimes_{\Lambda}\Z^k)$
satisfying the following:

Suppose that there is a unitary $v\in U(A)$ and a \morp\, $L: (C\rtimes_{\Lambda_{k-1}}\Z^{k-1})\otimes C(\T)\to A$
such
that
\beq\label{dp3-1}
&&\|[h(a), \, v]\|<\dt\tforal a\in {\cal  G},
{\rm{bott}}_1(h, v)|_{\cal  P}=0\andeqn\\\label{dp3-1+}
&&\tau\circ h|_{(C\rtimes_{\Lambda_{k-1}}\Z^{k-1})\otimes C(\T)}\approx_{\eta}\tau \circ L\,\,\,\text{on}\,\,\,{\bar {\cal G}},
\eneq
and
\beq\label{dp3-2}
L\approx_{\dt} h\,\,\,\text{on}\,\,\, {\cal G}\andeqn L(1\otimes z)\approx_{\dt} v
\eneq
and where $z\in C(\T)$ is the identity function on the unit
circle. Then there exists a unitary $W\in B\otimes{\cal U}^{k}$
such that
\beq\label{dp3-3}
\|[h^{(k)}(c), \,W]\|&<&\ep\rforal c\in {\cal F}\andeqn\\
W^*(v\otimes 1_{{\cal U}^{k}})W&\approx_{\ep}& h(1\otimes
z)\otimes 1_{{\cal U}^{k}},
\eneq
where $A\subset B\cong B\otimes{\cal U}$ is a unital simple
AF-algebra with a unique tracial state, with
$K_0(B)=\rho_B(K_0(B)$ which is also a subring of $\R$ and with
$1_A=1_B.$

\vspace{0.2in}

(4) We say that $C$ has the property (P4) with integer $k,$ if the
following holds:

Let $\Lambda: \Z^k\to Aut(C)$ be a \hm\, and let $A\cong A\otimes
{\cal U}$ be a unital simple AF-algebra with unique tracial state
$\tau$ and $\rho_A(K_0(A))=K_0(A).$ There exists at least one
unital monomorphism from $C\rtimes_{\Lambda}\Z^k$ into $A.$
For any $\ep>0$ and a finite subset ${\cal  F}\subset
C\rtimes_{\Lambda} \Z^k,$ there exists $\dt>0,$ a finite subset
${\cal G}\subset C\rtimes_\af \Z$ and a finite subset ${\cal
P}\subset \underline{K}(C\rtimes_{\Lambda}\Z^k)$ satisfying the
following:

 Suppose that $\phi: C\rtimes_{\Lambda}\Z^k\to A$ is a unital
monomorphism and suppose that there is a unitary $v\in U(A)$ such
that
\beq\label{dP4-1}
\|[\phi(a), \, v]\|<\dt\tforal a\in {\cal  G}\andeqn
{\rm{bott}}_1(\phi, v)|_{\cal  P}=0.
\eneq
Then there exists a continuous path of unitaries $\{v_t: t\in [0,
1]\}\subset U(B\otimes {\cal U}^{k} )$ such that
\beq\label{dP4-2}
v_0=v\otimes 1_{{\cal U}^{k}  },\,\,\, v_1=1_{A\otimes {\cal
U}^{k} },\,\,\,\|[\phi^{(k)}(a), v_t]\|<\ep
\eneq
for all $a\in {\cal  F}$ and $t\in [0,1],$ and
\beq\label{dP4-3}
{\rm{Length}}(\{v_t\})\le 2\pi+\min\{\ep,1\}.,
\eneq
where $A\subset B\cong B\otimes{\cal U}$ is a unital simple
AF-algebra with a unique tracial state, with
$K_0(B)=\rho_B(K_0(B)$ which is also a subring of $\R$ and with
$1_A=1_B.$

}

\end{df}

\begin{lem}\label{n2burg}
Let $C_0$ be a unital separable amenable \CA\, in ${\cal N}$ which
 has the
property {\rm (P4)} with integer $k-1>0$ and let $\Lambda: Z^k\to
Aut(C_0)$ be a \hm\,  such that  $\{\af_1,\af_2,...,\af_k\}$ forms
a set of generators of $\Lambda(\Z^k).$ Let $\Lambda_{k-1}:
\Z^{k-1}\to Aut(C_0)$ be defined by
$\{\af_1,\af_2,...,\af_{k-1}\}.$

For any $\ep>0$ and finite subset ${\cal F}\subset C_0\rtimes_{\Lambda_{k-1}}\Z^{k-1},$ there is
$\dt>0,$ a finite subset ${\cal G}\subset C_0\rtimes_{\Lambda_{k-1}}\Z^{k-1}$ and a finite subset
${\cal P}\subset K_1(C_0\rtimes_{\Lambda_{k-1}}\Z^{k-1})$
satisfying the following:

Suppose that $A\cong A\otimes {\cal U}$  is a  unital simple
AF-algebra with a unique tracial state and with
$K_0(A)=\rho_A(K_0(A)),$  suppose that $\phi:
C_0\rtimes_{\Lambda}\Z^k\to A$ is a unital monomorphism and
suppose that there is a unitary $v\in U(A)$ such that
\beq\label{n2bl-1}
\|[v,\, \phi(a)]\|<\dt,\,\,\,\tforal a\in {\cal  G} \tand
{\rm bott}_1(\phi, \, v)|_{\cal  P}=0
\eneq
then there is a unitary $w\in B\otimes {\cal U}^{k}$
such that
\beq\label{n2bl-2}
w^*(\phi^{(k-1)}|_{C_0\rtimes_{\Lambda_{k-1}}\Z^{k-1}}(a)\otimes
1_{{\cal U}^k})w&\approx_{\ep}&
\phi^{(k-1)}|_{C_0\rtimes_{\Lambda_{k-1}}\Z^{k-1}}(a)\otimes
1_{{\cal U}^{k}}\\
\hspace{-0.8in}\tforal a\in {\cal F} \tand\\
 w^*u(v\otimes 1_{{\cal U}^{k}})w&\approx_{\ep}& u,
\eneq
where $u=\phi^{(k)}(u_\af)$ and where $A\subset B\cong B\otimes
{\cal U}$ is a unital simple AF-algebra with a unique tracial
state, with $K_0(B)=\rho_B(K_0(B)$ which is also a subring of $\R$
and with $1_A=1_B.$
\end{lem}

\begin{proof}
Let $C=C_0\rtimes_{\Lambda_{k-1}}\Z^{k-1}.$ The proof is exactly
the same as that of \ref{2burg}. However, instead of using
$\phi',$ we use $\phi^{(k-1)}$ and instead of applying the Basic
Homotopy Lemma for AH-algebras, we apply the property (P4) with
integer $k-1.$

\end{proof}

\begin{thm}\label{n2Tuni}
Let $C_0$ be a unital separable amenable \CA\, in ${\cal N}$ which
has the property {\rm (P1)} with integer $k-1$ and {\rm (P4)} with
integer $k-1>0.$

Let $\Lambda: \Z^k\to Aut(C)$ be a \hm\,  and let $A\cong A\otimes
{\cal U}$ be a unital simple AF-algebra with a unique tracial
state $\tau$ and $K_0(A)=\rho_A(K_0(A)).$ Suppose that $\{\af_1,\,
\af_2,...,\af_k\}$ forms a set of generators for $\Lambda(\Z^k).$
Denote $\Lambda_{k-1}: \Z^{k-1}\to Aut(C_0)$ be defined by
$\{\af_1,\af_2,...,\af_{k-1}\}.$

 Suppose that
$\phi_1, \phi_2: C\rtimes_{\Lambda}\Z^k\to A$ are two unital
monomorphisms such that
\beq\label{n2Tu-1}
\tau\circ \phi_1=\tau\circ \phi_2.
\eneq
Suppose also that
\beq\label{n2Tu-1+}
R_{\phi_1\circ j_0, \phi_2\circ
j_0}(K_1(C_0\rtimes_{\Lambda_{k-1}}\Z^{k-1}))\subset
\rho_A(K_0(A)),
\eneq
where $j_0: C_0\rtimes_{\Lambda_{k-1}}\Z^{k-1}\to
C_0\rtimes_{\Lambda} \Z^k$ is the embedding.

  Then
there exists a sequence of unitaries $\{w_n\}\subset U(B\otimes
{\cal U}^{k})$ such that
\beq\label{n2Tu-2}
\lim_{n\to\infty}{\rm ad}\, w_n\circ
\psi_1^{(k)}(a)=\psi_2^{(k)}(a)\tforal a\in C\rtimes_{\af}\Z,
\eneq
where $A\subset B\cong B\otimes{\cal U}$ is a unital simple
AF-algebra with a unique tracial state, with
$K_0(B)=\rho_B(K_0(B)$ which is also a subring of $\R$ and with
$1_A=1_B.$
\end{thm}

\begin{proof}
Let $C=C_0\rtimes_{\Lambda_{k-1}}\Z^{k-1}.$ The proof of this is
exactly the same as that of \ref{2Tuni}. However, instead of
applying 10.7 of \cite{Lnasy}, (where AH-algebra is used) we use
the property (P1), and instead of applying \ref{2burg}, we apply \ref{n2burg}.

\end{proof}

\begin{thm}\label{nbott}
Let $C_{00}$ be a unital separable amenable \CA\, in ${\cal N}$
which has the property {\rm (P1)} with integer $k-1.$

 Then $C_{00}$ has the property
{\rm (P2)} with $k.$
\end{thm}

\begin{proof}
Let $\Lambda:\Z^k\to Aut(C_{00})$ be \hm, let
$\{\af_1,\af_2,...,\af_k\}$ forms a set of generators for
$\Lambda(\Z^k)$ and let $\Lambda_{k-1}: \Z^{k-1}\to Aut(C_{00})$
defined by $\{\af_1,\af_2,...,\af_{k-1}\}.$  Define
$C_0=C_{00}\rtimes_{\Lambda_{k-1}}\Z^{k-1}.$  The proof is very
much the same of that \ref{Bott1} but we apply the properties (P1)
instead of using the assumption that $C_0$ is an AH-algebra. We
proceed the same proof until (\ref{Bot-7}).  Note that, since
$C_{00}$ satisfies property (P1) with integer $k-1,$ \ref{n2Tuni}
holds. Next, instead of applying \ref{2Tuni}, we apply
\ref{n2Tuni} to obtain a unital simple AF-algebra $B_1\supset B$
($B$ as in the proof of \ref{Bott1}) with a unique tracial state,
with $K_0(B_1)=\rho_B(K_0(B_1))$ which is also a subring of $\R$
and with $1_B=1_{B_1}$ for which $B_1\cong B_1\otimes {\cal U},$
and a unitary $w\in B_1\otimes {\cal U}^k$ such that
(\ref{Bot-8})holds with $\psi_1$ replaced $(h\circ j_1)^{(k)}$ and
$\psi$ replaced by $j^{(k)}.$ Then (\ref{Bot-9}) follows (with
$\psi$ replaced by $j^{(k)}$).

\end{proof}

\begin{lem}\label{FKey}
Let $C_{00}$ be a unital separable amenable \CA\, in ${\cal N}$
which has the property {\rm (P2)} with $k-1,$  let $\Lambda:
\Z^k\to Aut(C_{00})$ be a \hm\, such that
$\{\af_1,\af_2,...,\af_k\}$ forms a set of generators of
$\Lambda(\Z^k)$ and let $\Lambda_{k-1}:\Z^{k-1}\to Aut(C_{00})$ be
defined $\{\af_1,\af_2,...,\af_{k-1}\}.$

Let $C_0=C_{00}\rtimes_{\Lambda_{k-1}}\Z^{k-1},$ let $A\cong
A\otimes {\cal U}$ be a unital simple AF-algebra with a unique
tracial state $\tau$ and with $K_0(A)=\rho_A(K_0(A))$ and let $j:
C_{00}\rtimes_{\Lambda}\Z^k\to A$ be a unital monomorphism, let
$\ep>0,$ ${\cal F}\subset C_0$ be a finite subset and ${\cal
P}\subset \underline{K}(C_0)$ be a finite subset.
 There
is $\dt>0,$ a finite subset ${\cal  G}\subset C_0,$ and integer
$K$ and a finite subset ${\cal  Q}\subset
\underline{K}(C_{00}\rtimes_{\Lambda}\Z^k)$ satisfying the
following:

Suppose that $v\in U(A)$ is a unitary such that
\beq\label{Fkey1}
\|[j(a),\, v]\|<\dt\tforal a\in {\cal  G}\andeqn
\eneq
\beq\label{Fkey2}
[L]|_{\cal  Q}=[j]|_{\cal  Q},
\eneq
where $L: C_{00}\rtimes_{\Lambda}\Z^k\to A$ is a \morp\, such that
\beq\label{Fkey3}
\|L(\sum_{k=-K}^{K}
f_ku_\af^k)-\sum_{k=-K}^{K}j(f_k)(j(u_\af)v)^k\|<\dt
\eneq
for all $f_k\in {\cal G}.$
Then there exists a unitary $w\in U(B\otimes {\cal U}^{k})$ such
that
\beq\label{Fkey4}
&&\|[j^{(k)}(a),\, w]\|<\ep\tforal a\in {\cal F}\andeqn\\
&&{\rm{bott}}_1(j^{(k)}|_{C_0},\, j^{(k)}(u_\af)^*w^*
j^{(k)}(u_\af)(v\otimes 1_{{\cal U}^{k}})w)|_{\cal  P}=0,
\eneq
 where
$A\subset B\cong B\otimes {\cal U}$ is a unital simple AF-algebra
with a unique tracial state, with $K_0(B)=\rho_B(K_0(B))$ which is
also a subring of $\R$ and with $1_A=1_B.$

\end{lem}

\begin{proof}
The proof of this lemma is almost the same as that of \ref{Key}
but we apply the property (P2) (with $k-1$) instead of Lemma 7.3
of \cite{Lnasy}. We proceed as follows.

We will keep the first paragraph of that proof  and keep the line
which defines ${\cal P}.$
 Let $\eta>0$ (in place of $\dt$),
$\{y_1,y_2,...,y_m\}\subset K_1(C_0))$ be required by (P2) for
$\ep_0/2>0,$  $\{x_1,x_2,...,x_k\}$ and ${\cal F}_1.$ We then
remove the paragraph  related to Lemma 7.3 of \cite{Lnasy} right
after the line which defines ${\cal P}.$

We then continue the same proof until the line right after (\ref{key-14}).
We replace the next line
by the following:

 By (\ref{key-4}) and by
property (P2) with $k-1$ for $C_{00},$ there is a unitary $w\in
B\otimes {\cal U}^{k}$ such that (\ref{key-14+}) holds with $j$
replaced by $j^{(k)},$ where $A\subset B\cong B\otimes{\cal U}$ is
a unital simple AF-algebra with a unique tracial state, with
$K_0(B)=\rho_B(K_0(B)$ which is also a subring of $\R$ and with
$1_A=1_B.$
 The rest of the proof remains the same by replacing
$j$ by $j^{(k)}$ and $v$ by $v\otimes 1_{{\cal U}^k}.$

\end{proof}

\begin{thm}\label{TP3}
Let $C_{0}$ be a unital separable amenable \CA\, in ${\cal N}$
 which
has property {\rm (P1)}, {\rm (P2)}, {\rm (P3)} and {\rm (P4)}
with $k-1.$

Then $C_{0}$ has the property {\rm (P3)} with integer $k$ .

\end{thm}

\begin{proof}

The proof is the same as that of \ref{S1uni}. The main difference
is that  we will use property (P3) instead of \ref{2QUni}.

We proceed as follows:

Let $\Lambda: \Z^k\to Aut(C_{0})$ be a monomorphism such that
 $\{\af_1,\af_2,...,\af_k\}$ forms a set of generators for
 $\Lambda(\Z^k)$ and let $\Lambda_{j}:\Z^{j}\to Aut(C_{0})$
 defined by $\{\af_1, \af_2,...,\af_{j}\}$ ($1\le j\le k-1$).

 Define $C=C_0\rtimes_{\Lambda_{k-1}}\Z^{k-1}.$

Since $C_0$ has property (P1), (P2), (P3) and (P4) (with $k-1$),
\ref{n2burg} and \ref{FKey} hold.

We will not make changes from the beginning of the proof to
(\ref{lh-2-1+}) and keep the definition of $\tau_0$ with only a
couple of exceptions: We will not apply \ref{2burg}, but apply
\ref{n2burg}, not \ref{Key} but \ref{FKey}. So we replace
\ref{2burg} by \ref{n2burg} and replace \ref{Key} by \ref{FKey} in
these lines. We also need, however, replace $\af$ by $\af_k.$

We may also assume that ${\cal G}_2={\cal G}_0\cup\{1\otimes z\}$
for some finite subset ${\cal G}_0\subset C$ and ${\cal
G}_2'={\cal G}_0'\cup \{1\otimes z\}$ for some finite subset ${\cal G}_0'.$

After we define $\tau_0$ (right after (\ref{lh-2-1+})), we will
apply property (P3) with integer $k-1.$

Let $0<\dt_3<\dt_2$ (in place of $\dt$),  $\eta>0,$ ${\tilde {\cal
G}} \subset C_0\rtimes_{\Lambda_{k-1}}\Z^{k-1}$ (in place of
${\cal G}$) be a finite subset, ${\bar {\cal G}}_1\subset
(C_0\rtimes_{\Lambda_{k-2}}\Z^{k-2})\otimes C(\T)$ (in place of
${\bar {\cal G}}$) be a finite subset required and ${\cal
P}_4\subset K_1(C)$ be a finite subset required by (P3) (with
$k-1$) for $\min\{\dt_2/4,\ep/4\}$ and ${\cal G}_0'.$

As in the proof of \ref{S1uni}, without loss of generality, we may assume that ${\bar {\cal  G}}_1={\cal G}_3'\otimes {\cal G}_3'',$ where $\{1_C\}\subset {\cal G}_3'\subset C$ and $\{1_{C(\T)}, z\}\subset {\cal G}_3''\subset C(\T)$ are finite subsets.

We delete the lines after the definition of $\tau_0$ until the
definition of $j_0$ of the proof of \ref{S1uni}. Keep the
definition of $j_0$ and put ${\cal P}=(j_0)_{*1}({\cal P}_4)$ and
keep the same $\eta_1.$ We keep the next few lines but stop right
after (\ref{lh-4}) and replace (\ref{S1u1}) by (\ref{dp3-1}),
(\ref{S1u1+}) by (\ref{dp3-1+}) and (\ref{S1u2}) by (\ref{dp3-2}).
Since we assume that (\ref{dp3-2}) holds, we have
\beq\label{nS1-10}
\text{bott}_1(h|_{C} ,\, v)|_{{\cal P}_4}=0.
\eneq
By applying property (P3) (with $k-1$), we obtain a unitary
$w_0'\in B_1\otimes {\cal U}^{k-1}$ such that
\beq\label{nS1-11}
\|[(h|_C)^{(k-1)}(c),\, w_0']\|&<&\min\{\dt_2/4, \ep/2\}\rforal c\in {\cal G}_0'\andeqn \\
(w_0')^*(v\otimes 1_{{\cal U}^{k-1}})w_0'&\approx&
\hspace{-0.1in}_{\min\{\dt_2/4, \ep/2\}} h(1_{C_0}\otimes
z)\otimes 1_{{\cal U}^{k-1}},
\eneq
where $A\subset B_1\cong B_1\otimes {\cal U}$ is a unital simple
AF-algebra with a unique tracial state, with
$K_0(B_1)=\rho_{B_1}(K_0(B_1))$ which is also a subring of $\R$
and with $1_A=1_{B_1}.$

 Define $w_0=w_0'\otimes 1_{\cal U}.$ Then
\beq\label{nS1-12}
\|[h^{(k)}(c),\, w_0]\|&<&\min\{\dt_2/4, \ep/2\}\rforal c\in {\cal G}_0'\andeqn \\\label{nS1-12+}
(w_0)^*(v\otimes 1_{{\cal U}^{k}})w_0&\approx&\hspace{-0.1in}_{\min\{\dt_2/4, \ep/2\}}  h(1_{C_0}\otimes z)\otimes 1_{{\cal U}^{k}}
\eneq

Define $V=(h^{(k)}(u_{\af_k}))^*w_0^*h^{(k)}(u_{\af_k})w_0.$
We then continue the same proof as in \ref{S1uni} as follows.
We have
\beq\label{nS-13}
\|[h^{(k)}(a),\, V]\|<\min\{\ep/2, \dt_2/4\}\rforal a\in {\cal G}_2'.
\eneq
Thus, there exists a unital \morp\, $L_1: C\rtimes_{\Lambda}\Z^k\otimes C(\T)\to A\otimes {\cal U}^k$
which satisfies (\ref{lh-2-1}) (replacing $\af$ by $\af_k$ and $h$ by $h^{(k)}$) with $v'=V.$
Therefore
\beq\label{nS1-14}
L_1\approx_{\dt_2} h^{(k)}\,\,\,\text{on}\,\,\, {\cal G}_2''\andeqn L_1(u_{\af_k})\approx_{\dt_2}h^{(k)}(u_{\af_k})V.
\eneq
Note that
\beq\label{nS1-15}
[L_1]|_{{\cal P}_2}=[h^{(k)}]|_{{\cal P}_2}
\eneq

By the assumption (\ref{dp3-1}) and by (\ref{nS1-12+}),
\beq\label{nS1-16}
\text{bott}_1({\rm ad}\, w_0\circ h^{(k)},\, v\otimes 1_{{\cal U}^k}]\||_{{\cal P}_3} &=&
\text{bott}_1({\rm ad}\, w_0\circ h^{(k)},\, w_0^*(v\otimes 1_{{\cal U}^k})w_0)|_{{\cal P}_3}\\
&=&\text{bott}_1(h^{(k)},\, v)|_{{\cal P}_3}=0.
\eneq
It follows that (since ${\rm ker}\rho_{B_1}=\{0\}$)
\beq\label{nS1-17}
[L_1]|_{{\boldsymbol\bt}({\cal P}_3)}=0=[h^{(k)}]|_{{\boldsymbol\bt}({\cal P}_3)}.
\eneq

Therefore, by applying \ref{FKey}, there is a unitary $w_1\in
U((B_2\otimes {\cal U}^k)$ such that
\beq\label{nS1-18}
\|([h^{(k)})^{(k)}(a),\, w_1]\|&<&\dt_1/2\rforal a\in {\cal G}_1\andeqn\\
 \text{bott}_1((h^{(k)})^{(k)},\, V_1)|_{{\cal P}_1}=0,
 \eneq
 where
 $V_1=(h^{(k)})^{(k)}(u_{\af_k})^*w_1^*(h^{(k)})^{(k)}(u_{\af_k})(V\otimes 1_{{\cal
 U}^k})w_1$
 and where $B_1\otimes {\cal U}^k\subset B_2\cong B_2\otimes {\cal U}$ is a unital
 simple AF-algebra with a unique tracial state, with
 $K_0(B_2)=\rho_{B_2}(K_0(B_2))$ which is also a subring of $\R$
 and with $1_{B_1\otimes {\cal U}^k}=1_{B_2}.$
 Note that
 $$(h^{(k)})^{(k)}(u_{\af_k})(V\otimes 1_{{\cal U}^k})=w_0^*h^{(k)}(u_{\af_k})w_0\otimes 1_{{\cal U}^{k-1}}\otimes\imath(u_\sigma).$$
 It follows from \ref{n2burg} that there exists a unitary
 $w_2\in U(B_2\otimes {\cal U}^{2k})$ such that
 \beq\label{nS1-19}
 &&w_2^*(((h^{(k)})^{(k)})^{(k)}(c))w_2\approx_{\ep/2} ((h^{(k)})^{(k)})^{(k)}(a)\rforal a\in {\cal F}\\\label{nS1-20-}
 &&\hspace{-0.8in}w_2^*(w_1^*(w_0^*h^{(k)}(u_{\af_k})w_0\otimes 1_{{\cal U}^k})w_1)\otimes 1_{{\cal U}^k} \otimes\imath(u_{\af_k})\approx_{\ep/2}
 (h^{(k)})^{(k)}(u_{\af_k})\otimes 1_{{\cal U}^{k-1}}\otimes \imath(u_{\af_k}).
\eneq
Define $W_1=((w_0\otimes 1_{\cal U})w_1)\otimes 1_{\cal U})w_2.$
Then, by (\ref{nS1-12}) and (\ref{nS1-12+}),
\beq\label{nS1-20}
&&\|[((h^{(k)})^{(k)})^{(k)}(c),\, W_1]\|<\ep\rforal c\in {\cal
F}\andeqn\\
&& W_1^*(v\otimes 1_{{\cal U}^{3k}})W \approx_{\ep}
h(1\otimes z)\otimes 1_{{\cal U}^{3k}}.
\eneq

Finally, we apply the absorption lemma \ref{Abs} twice and
\ref{3abs} once.

\end{proof}

\begin{lem}\label{nlhom}
Let $C$ be a unital separable amenable \CA\, in ${\cal N}$ and
 let $\Lambda:
\Z^k\to Aut(C)$ be a \hm. Suppose that $C$ has property (P3) for
$k.$ Then $C$ has the property {\rm (P4)} with $k.$
\end{lem}

\begin{proof}
The proof follows from that of \ref{Lhom} almost verbatim. We
replace the application of \ref{S1uni} by the application of
property (P3). So in (\ref{nhl-13}), $W$ should be in $U(B_1
\otimes {\cal U}^k),$ ${\bar v}$ should be replaced by ${\bar
v}\otimes 1_{{\cal U}^k},$ $h(1\otimes z)\otimes 1_{\cal U}$
replaced by $h(1\otimes z)\otimes 1_{{\cal U}^k}$ and $\psi$
should be replaced by $\phi^{(k)},$ where $A\otimes {\cal
U}\subset B\cong B\otimes {\cal U}$ is a unital simple AF-algebra
with a unique tracial state, with $K_0(B)=\rho_B(K_0(B))$ which is
also a subring of $\R$ and with $1_{A\otimes {\cal U}}=1_B.$ For
the rest of the proof, we will continue use this $W$ and replace
$1_{\cal U}$ by $1_{{\cal U}^k},$ $\psi$ by $\phi^{(k)}$ and $V$
by $V\otimes 1_{{\cal U}^k}$ (and $\lambda(t)$ and $v_t$ are in
$B\otimes {\cal U}^k$).

\end{proof}

\begin{thm}\label{Lind}
Let  $C$ be  a unital separable amenable \CA\, in ${\cal N}$ and
let $\Lambda: \Z^k\to Aut(C)$ be a monomorphism.

Suppose that $C$ satisfies property {\rm (P1)}, {\rm (P2)}, {\rm
(P3)} and {\rm (P4)} for $k-1.$ Then $C$ has property {\rm (P1)},
{\rm (P2)}, {\rm (P3)} and {\rm (P4)} for $k.$

\end{thm}

\begin{proof}
Since $C$ has property (P1) (with integer $k-1$), \ref{n2Tuni}
holds. It follows from \ref{nbott} and \ref{TP3} that $C$ has the
property (P2) and (P3) for $k.$ By \ref{nlhom}, $C$ also has the
property (P4) for $k.$ It remains to show that $C$ has property
(P1) for $k.$

The proof is exactly the same as that of \ref{TAsy}.

Suppose that $\{\af_1,\af_2,...,\af_k\}$ forms a set of generators
for $\Lambda(\Z^k).$ Let $\Lambda_{k-1}: \Z^{k-1}\to Aut(C)$ be
defined by $\{\af_1, \af_2,...,\af_{k-1}.$

In the proof of \ref{TAsy}, we need to replace $C$ by
$C\rtimes_{\Lambda_{k-1}}\Z^{k-1},$ $\af$ by $\af_k,$ $\psi_1$ and
$\psi_2$ by $\phi_1^{(k)}$ and $\phi_2^{(k)},$ respectively. In
particular, (\ref{asy-4}) holds for $\phi_1^{(k)}$ and
$\phi_2^{(k)}$ (instead of $\psi_1$ and $\psi_2$). We also replace
$A\otimes {\cal U}$ by $B_1\otimes {\cal U}^k$ for some unital
simple AF-algebra $B_1$ which has the properties: (1) $A\subset
B_1\cong B_1\otimes {\cal U},$ (2) $K_0(B_1)=\rho_{B_1}(K_0(B_1))$
is a subring of $\R,$ (3) $1_{A}=1_{B_1}.$


We will prepare to apply (P4) (with $k$) instead of \ref{Lhom}. We
will apply property (P2) (with $k$). Hence (\ref{asy-23}) holds by
replacing $\psi_2'$ by $(\phi_2^{(k)})^{(k)}.$ We proceed the same
calculation and apply property (P4) (with $k$) to obtain
(\ref{asy-32}) with $\psi_2''$ being  replaced by
$((\phi_2^{(k)})^{(k)})^{(k)}$ and with $v_n(t)\in B_2\otimes
{\cal U}^{k},$ for some unital simple AF-algebra $B_2\cong
B_2\otimes {\cal U}$ which satisfies: (1) $B_1\otimes {\cal
U}^k\subset B_2,$ (2) $K_0(B_2)=\rho_{B_2}(K_0(B_2))$ which is a
subring of $\R,$ (3) $1_{B_1\otimes {\cal U}^k}=1_{B_2}.$

We then obtain (\ref{asy-33}) with $\psi_1''$ and $\psi_2''$ being
replaced by $(\phi_1^{(k)})^{(k)})^{(k)}$ and
$(\phi_2^{(k)})^{(k)})^{(k)},$ and with $u(t)\in B_2\otimes {\cal
U}^{2k}.$

Finally, to obtain (\ref{Asy2}), we apply the absorption lemma
\ref{Abs} twice.

\end{proof}

\begin{thm}\label{asy-emb}
Let $C$ be a unital separable amenable \CA\,  in ${\cal N}.$
Suppose that $\Lambda: \Z^k\to Aut(C)$ is a \hm\, so that
$\{\af_1,\af_2,...,\af_k\}$ forms a set of generators for
$\Lambda(\Z^k)$ and suppose that $\Lambda_{k-1}: \Z^{k-1}\to
Aut(C)$ defined by $\{\af_1, \af_2,...,\af_{k-1}\}.$

Suppose that there is a unital monomorphism $\phi:
C\rtimes_{\Lambda_{k-1}}\Z^{k-1}\to A,$ where $A\cong A\otimes
{\cal U}$ is a unital simple AF-algebra with a unique tracial
state $\tau$ and with $K_0(A)=\rho_A(K_0(A))$ such that $\tau\circ
\phi$ is $\af_k$-invariant.

Suppose $C$ has property {\rm (P1)} with $k-1.$  Then there is a
unital monomorphism $\psi: C\rtimes_{\Lambda}\Z^k\to B$ for some
unital simple AF-algebra with a unique tracial state $\tau_1$ and
with $K_0(B)=\rho_B(K_0(B))$ such that
\beq\label{asyemb1}
\tau_1\circ \psi\circ j_0=\tau\circ \phi,
\eneq
where $j_0: C\rtimes_{\Lambda_{k-1}}\Z^{k-1}\to
C\rtimes_{\Lambda}\Z^k$ is the embedding.
\end{thm}

\begin{proof}
The proof is an application of  Theorem 5.5 of \cite{Lnasy} exactly as in the proof of \ref{1MT}
(see also \ref{Tz2}).

\end{proof}

\begin{thm}\label{FMT}
Let $C$ be a unital separable AH-algebra and let $\Lambda: \Z^k\to
Aut(C)$ be a monomorphism. Then $C$ can be embedded into a unital
simple AF-algebra if and only if $C$ admits a faithful
$\Lambda$-invariant tracial state.

\end{thm}

\begin{proof}
It suffices to show the ``if part" of the theorem.

Suppose that $\{\af_1, \af_2,...,\af_k\}$ forms a set of generators for $\Lambda(\Z^k).$
Denote by $\Lambda_j: Z^j\to Aut(C)$ defined by $\{\af_1, \af_2,...,\af_j\},$ $j=1,2,...,k.$

 By \ref{Bott1},  \ref{S1uni}, \ref{Lhom}  and \ref{TAsy}, $C$ has property
(P1) (P2), (P3), (P4) for $k=1.$ Theorem \ref{asy-emb} (or \ref{casy})
implies that there is a unital monomorphism $\phi_2: C\rtimes_{\Lambda_2}\Z^2\to A_2,$ where
$A_2\cong A_2\otimes {\cal U}$ is a unital simple AF-algebra with a unique tracial state $\tau_2$ and with $K_0(A_2)=\rho_{A_2}(K_0(A_2))$ such that
\beq\label{FMT-1-1}
\tau_2\circ \phi_2\circ j_1=t,
\eneq
where $t$ is the chosen faithful $\Lambda$-invariant tracial state on $C.$ In particular,
\beq\label{FMT-1-1+}
\tau_2\circ \phi_2=\tau_2\circ \phi_2\circ \af_i,\,\,\, i=1,2,...,k.
\eneq

We use the induction on $k$ to
prove the theorem.

Suppose that we have shown the following:

(1)  $C$ has the property (P1), (P2), (P3) and (P4) for $i-1.$

(2) There is a monomorphism $\phi_{i-1}: C\rtimes_{\Lambda_{i-1}}\
Z^{i-1}\to A_{i-1},$ where $A_{i-1}\cong A_{i-1}\otimes {\cal U}$
is a unital simple AF-algebra with a unique tracial state
$\tau_{i-1}$ and with $K_0(A_{i-1})=\rho_{A_{i-1}}(K_0(A_{i-1}))$
such that
\beq\label{FMT-1}
\tau_{i-1}\circ \phi_{i-1}=\tau_{i-1}\circ \phi\circ\af_l,\,\,\, l=1,2,...,k.
\eneq

It follows from \ref{Lind} that $C$ has property (P1), (P2),(P3)
and (P4) for $i.$ It follows from \ref{asy-emb} that there is a
unital monomorphism $\phi_{i}: C\rtimes_{\Lambda_i}\Z^i\to
A_{i}\cong A_{i}\otimes {\cal U}$ for some unital simple
AF-algebra $A_{i}$ with a unique tracial state $\tau_i$ and with
$K_0(A_{i})=\rho_{A_{i}}(K_0(A_{i}))$ such that
\beq\label{FMT-2}
\tau_i\circ \phi_i\circ j_i=\tau_{i-1}\circ \phi_{i-1},
\eneq
where $j_i: C\rtimes_{\Lambda_{i-1}}\Z^{i-1}\to C\rtimes_{\Lambda_i}\Z^i$
is the embedding.
By (\ref{FMT-1}), we have
\beq\label{FMT-3}
\tau_i\circ \phi_i=\tau_i\circ \phi_i\circ \af_l,\,\,\,
l=1,2,....,k.
\eneq

Thus, there is a unital monomorphism $\phi_k: C\rtimes_{\Lambda}\Z^k\to A_k,$ where
$A_k\cong A_k\otimes {\cal U}$ is a unital simple AF-algebra with a unique tracial state $\tau$ and with $K_0(A_k)=\rho_{A_k}(K_0(A_k)).$

\end{proof}

\begin{thm}\label{FFMT}
Let $C$ be a unital AH-algebra and let $G$ be a finitely generated
(discrete)abelian group. Suppose that $\Lambda: G\to Aut(C)$ is a
\hm. Then $C\rtimes_{\Lambda} G$ can be embedded into a unital
simple AF-algebra if and only if $C$ admits a faithful
$\Lambda$-invariant tracial state.
\end{thm}

\begin{proof}
We use an argument of Nate Brown (see Remark 11.10 of \cite{Bn2}).
Write $G=\Z^k\oplus G_0,$ where $G_0$ is a finite group.
By the Green's theorem (Cor. 2.8 of \cite{Gr}), there is a monomorphism $h: C\rtimes_{\Lambda} G\to
(C\rtimes_{\Lambda|_{\Z^k}}\Z^k)\otimes {\cal K}.$ Since $C\rtimes_{\Lambda|_{\Z^k}}\Z^k$ has a faithful
$\Lambda|_{\Z^k}$ -invariant tracial state, there is a unital monomorphism $\phi: C\rtimes_{\Lambda|_{\Z^k}}\Z^k\to A$ for some unital simple AF-algebra.
Then $\phi$ gives a monomorphism $\phi_1: (C\rtimes_{\Lambda|_{\Z^k}}\Z^k)\otimes {\cal K}\to A\otimes {\cal K}.$
Thus $\phi_1\circ h$ gives a monomorphism from $C\rtimes_{\Lambda} G$ into
$A\otimes {\cal K}.$ Let $e=\phi_1\circ h(1).$ Then $e$ is a projection.
Note that $e(A\otimes {\cal K})e$ is a unital simple AF-algebra.

\end{proof}

\begin{cor}\label{voicc}
Let $X$ be a compact metric space and $\Lambda: \Z^k\to Homeo(X)$
be a \hm, where $Homeo(X)$ is the group of all homeomorphisms on
$X.$ Then the crossed product $C(X)\rtimes_{\Lambda}\Z^k$ can be
embedded into a unital simple AF-algebra if and only if there is a
strictly positive $\Lambda$-invariant probability Borel measure.
\end{cor}

\end{document}